\documentclass{IEEEtran}
\usepackage{amsmath,amssymb,amsfonts}
\usepackage{cuted}
\usepackage{graphicx}
\usepackage{booktabs}
\usepackage[usenames,dvipsnames]{xcolor}
\usepackage{soul}
\usepackage{balance}
\usepackage{amsthm}

\newcommand{\1}[0]{\boldsymbol{1}}
\newcommand{\E}[0]{\mathbf E}
\newcommand{\R}[0]{\mathbf R}
\newcommand{\N}[0]{\mathbf N}
\renewcommand{\P}[0]{\mathcal P}
\newcommand{\D}[0]{\mathcal D}
\newcommand{\Z}[0]{\mathbf Z}

\def\conv{\mathrm{conv}}

\newcommand{\vv}[1]{\boldsymbol{#1}}

\newcommand{\iter}[2]{#1_{#2}}

\newtheorem{theorem}{Theorem}
\newtheorem{lemma}{Lemma}
\newtheorem{corollary}{Corollary}
\newtheorem{assumption}{Assumption}

\newtheorem{definition}{Definition}

\graphicspath{{Figures/}}

\title{A Convex Optimization Approach to \\ Discrete Optimal Control}
\author{V\'ictor Valls and Douglas J. Leith\\ Trinity College Dublin \thanks{This work was supported by Science Foundation Ireland under Grant No. 11/PI/1177.}}

\begin{document}

\maketitle

\begin{abstract}
In this paper, we bring the celebrated max-weight features (myopic and discrete actions) to mainstream convex optimization. Myopic actions are important in control because decisions need to be made in an online manner and without knowledge of future events, and discrete actions because many systems have a finite (so non-convex) number of control decisions. For example, whether to transmit a packet or not in communication networks.
Our results show that these two features can be encompassed in the subgradient method for the Lagrange dual problem by the use of stochastic and $\epsilon$-subgradients. 
One of the appealing features of our approach is that it decouples the choice of a control action from a specific choice of subgradient, which allows us to design control policies without changing the underlying convex updates. 
Two classes of discrete control policies are presented: one that can make discrete actions by looking only at the system's current state, and another that selects actions using blocks. The latter class is useful for handling systems that have constraints on the order in which actions are selected. 
\end{abstract}

\begin{IEEEkeywords}
approximate optimization, convex optimization, discrete and stochastic control, max-weight scheduling, subgradient methods.
\end{IEEEkeywords}

\section{Introduction}

\IEEEPARstart{C}{onvexity} plays a central role in mathematical optimization from both a theoretical and practical point of view. 
Some of the advantages of formulating a problem as a convex optimization are that there exist numerical methods that can solve the optimization problem in a reliable and efficient manner, and that a solution is always a global optimum. When an optimization problem is not of the convex type, then one enters into the realm of non-convex optimization where the specific structure of a problem must be exploited to obtain a solution, often not necessarily the optimal one.

In some special cases there exist algorithms that can find optimal solutions to non-convex problems. One is the case of max-weight scheduling: an algorithm initially devised for scheduling packets in queueing networks which has received much attention in the networking and control communities in recent years. In short, max-weight was proposed by Tassiulas and Ephremides in their seminal paper \cite{TE92}. It considers a network of interconnected queues in a slotted time system where packets arrive in each time slot and a network controller has to make a discrete (so non-convex) scheduling decision as to which packets to serve from each of the queues.
Appealing features of max-weight are that the action set matches the actual decision variables (namely, do we transmit or not); that scheduling decisions are made without previous knowledge of the mean packet arrival rate in the system (myopic actions); and that it can stabilize the system (maximize the throughput) whenever there exists a policy or algorithm that would do so. These features have made max-weight well-liked in the community and have fostered the design of extensions that consider (convex) utility functions, systems with time-varying channel capacity and connectivity, heavy-tailed traffic, \emph{etc.} Similarly, max-weight has been brought to other areas beyond communication networks including traffic signal control \cite{WUW12}, cloud computing \cite{MSY12}, economics \cite{Nee10c}, \emph{etc.}, and has become a key tool for discrete decision making in queueing systems.

However, there are some downsides. The success of max-weight has produced so many variants of the algorithm that the state of the art is becoming not only increasingly sophisticated but also increasingly complex. Furthermore, it is often not clear how to combine the different variants (\emph{e.g.}, utility function minimization with heavy-tailed traffic) since the design of a new control or scheduling policy usually involves employing a new proof or exploiting the special structure of a problem. There is a need for abstraction and to put concepts into a unified theoretical framework. 
While this has been attempted in previous works by means of establishing a connection between max-weight and dual subgradient methods in convex optimization, most of the works  \cite{LS04, LS06, ES06, ES07, CLCD06} have focused on specific congestion control applications and obtained discrete control actions as a result of exploiting the special structure of the problem. In particular, discrete actions are possible because the primal problem allows decomposition and a partial subgradient (or schedule) can be obtained as a result of minimizing a linear program.\footnote{A solution of a linear program always lies in an extreme point of a polytope, and that point is matched with a discrete (scheduling) decision.}
Further, the aforementioned works are deeply rooted in Lyapunov or fluid limit techniques and convergence/stability is only guaranteed asymptotically. 
That is, they do not make quantitative statements about the system state in finite time; for instance, provide upper and lower bounds on the optimality of the objective function, or a bound on the amount of constraint violation.
Because of all the above, the body of work on max-weight approaches is still largely separate from the \emph{mainstream} literature on convex optimization.

In this paper, we abstract the celebrated max-weight features (myopic and discrete actions) and make them available in standard convex optimization.
In particular, our approach consists of formulating the Lagrange dual problem and equipping the subgradient method with a perturbation scheme that can be regarded as using stochastic and $\epsilon$-subgradients. Stochastic subgradients are useful to capture the randomness in the system, while $\epsilon$-subgradients allow us to \emph{decouple} the choice of subgradient from the selection of a control action, \emph{i.e.}, they provide us with flexibility as to how to select control actions. The latter is important because (i) we can design a new control policy without having to prove again the convergence of the whole algorithm, and (ii) it eases the design of policies that model the characteristics of more complex systems.
For example, policies that select actions from a finite set or that have constraints associated with selecting certain subsets of actions. The finiteness of the action set is of particular significance from a theoretical point of view because we are allowing convex optimization to make non-convex updates, and from a practical point of view because many systems, such as computers, make decisions in a discrete-like manner.  

The main contributions of the paper are summarized in the following:

\begin{itemize}  
\setlength\itemsep{0.5em}

\item [(i)] \emph{Unifying Framework:}  Our analysis brings the celebrated max-weight features to mainstream convex optimization. In particular, they can be encompassed in the subgradient method for the Lagrange dual problem by using $\delta_k$ and $\epsilon_k$ perturbations. The analysis is presented in a general form and provides different types of convergence depending on the statistical properties of the perturbations, including bounds that are not asymptotic.

\item[(ii)] \emph{General Control Policies:} Our analysis clearly separates the selection of a subgradient from a particular choice of control action and establishes the fundamental properties that a control policy should satisfy for it to be optimal. 

\item [(iii)] \emph{Discrete Control Policies:} We develop two classes of control policies that allow us to use action sets with a finite number of actions (\emph{i.e.}, the action sets are not convex). One that is able to make discrete actions by looking only at the system's current state, and another that selects actions using blocks. The latter class is useful for handling systems that have constraints on how actions can be selected.
\end{itemize}

The rest of the paper is organized as follows. We start with the preliminaries, which cover the notation and some background material. In Section \ref{sec:main}, we study the convergence of the dual subgradient method under a ($\delta_k, \epsilon_k$) perturbation scheme, and in Section \ref{sec:actions} how the $\epsilon_k$ perturbations can be used to equip the dual subgradient method with discrete actions.
Section \ref{sec:remarks} provides some remarks and discussion, and Section \ref{sec:example}  illustrates the results with an example that considers discrete scheduling decisions with constraints. Finally, in Section \ref{sec:relatedwork} we provide an overview of the state of the art, and compare it with our work. All the proofs are in the appendices.

\section{Preliminaries}

We start by introducing the notation, the standard convex optimization problem setup, and the subgradient method for the Lagrange dual problem.

\subsection{Notation}
The sets of natural, integers and real numbers are denoted by $\N$, $\Z$ and $\R$. We use $\R_+$ and $\R^n$ to denote the set of nonnegative real numbers and $n$-dimensional real vectors. Similarly, we use $\R^{m \times n}$ to denote the set of $m \times n$ real matrices. Vectors and matrices are usually written, respectively, in lower and upper case, and all vectors are in column form.  The transpose of a vector $x \in \R^n$ is indicated with $x^T $, and we use $\1$ to indicate the all ones vector. The Euclidean, $\ell_1$ and $\ell_\infty$ norms of a vector $x \in \R^n$ are indicated, respectively, with $\|x\|_2$, $\| x \|_1$ and $\| x \|_\infty$.

Since we usually work with sequences we will use a subscript to indicate an element in a sequence, and parenthesis to indicate an element in a vector. For example, for a sequence $\{ x_k \}$ of vectors from $\R^n$ we have that $x_k = [x_k(1), \dots, x_k(n) ]^T $ where $x_k(j), \ j=1,\dots,n$ is the $j$'th component of the $k$'th vector in the sequence. For two vectors $x, y  \in \R^n$ we write $x \succ y$ when $x{(j)} > y{(j)}$ for all $j= 1,\dots,n$, and $x \succeq y$ when $x(j) \ge  y(j)$. We use $[\cdot]^{+}$ to denote the projection of a vector $x \in  \R^n$ onto the nonnegative orthant, \emph{i.e.}, $[x]^+ = [\max \{x(1), 0 \}, \dots, \max \{ x(n), 0\}]^T$.

\subsection{Convex Optimization Problem}
Consider the standard constrained convex optimization problem $\P$
\begin{align*}
\begin{array}{ll}
\underset{x \in X}{\text{minimize}} & f(x)\\
\mbox{\text{subject to}} &  g_j(x) \le 0 \qquad j=1,\dots,m
\end{array}
\end{align*}
where $f, g_j : X \to {\R }$ are convex functions and $X$ is a convex subset from $\R^n$. We will assume that set $X_0 := \{ x \in X \mid  g_j(x) \le 0, \ j=1,\dots,m\} \ne \emptyset$, and so problem $\P$ is feasible. Also, and using standard notation, we define $f^\star := \min_{x \in X_0} f(x)$ and $x^\star \in \arg \min_{x \in X_0} f(x)$. 

We can transform problem $\P$ into an unconstrained convex optimization by applying a (Lagrange) relaxation on the constraints. The Lagrange dual function associated to problem $\P$ is given by
\begin{align*}
h(\lambda)=\inf_{x \in X} L(x,\lambda) = \inf_{x \in X} \left\{ f(x) + \lambda^T  g(x) \right\},
\end{align*}
where $g(x) = [ g_1(x),\dots, g_m(x)]^T $, and $\lambda \in \R^m_+$ is a vector of Lagrange multipliers. 
Function $h$ is concave \cite[Chapter 5]{BV04} and so we can cast the following unconstrained\footnote{The Lagrange dual function is equal to $-\infty$ when $\lambda \prec 0$.} concave maximization problem $\D$
\begin{align*}
&\underset{\lambda \succeq 0}{\text{maximize }} \quad h(\lambda)
\end{align*}
where $h(\lambda^\star) = f^\star$ with $\lambda^\star \in \arg \max_{\lambda \succeq 0} h(\lambda)$ when strong duality holds. That is, solving problem $\P$ is equivalent to solving problem $\D$. 
A sufficient condition for strong duality to hold is the following.
\begin{assumption}[Slater condition] 
\label{th:slater_general}
$X_0$ is non-empty. There exists a vector $x \in X$ such that $g_j(x) < 0$ for all $j = 1,\dots,m$.
\end{assumption}
%

\subsection{Classic Subgradient Method}
Problem $\D$ can be solved using the subgradient method. One of the motivations for using this iterative method is that it allows the Lagrange dual function to be nondifferentiable, and so it imposes few requirements on the objective function and constraints.
Another motivation for using the subgradient method is that when the Lagrangian has favorable structure, then the algorithm can be implemented in a distributed manner, and therefore can be used to solve large-scale problems.
Nonetheless, in this work, the primary motivation for using the subgradient method in the Lagrange dual problem is that it allows us to handle perturbations on the constraints. As we will show in Section \ref{sec:main}, this will be key to handling resource allocation problems where the resources that need to be allocated are not known in advance.

\subsubsection{Iteration}In short, the subgradient method for the Lagrange dual problem consists of the following update:
\begin{align*}
\iter{\lambda}{k+1} =  [\iter{\lambda}{k} + \alpha \partial h (\lambda_k)]^+ \qquad k=1,2,\dots
\end{align*}
where $\iter{\lambda}{1} \in \R^m_+$, $\partial h (\lambda_k)$ is a subgradient in the subdifferential of $h$ at point $\iter{\lambda}{k}$ and $\alpha > 0$ a constant step size. The classic subgradient method can make use of more complex step sizes, but constant step size will play an important role in our analysis, and it is extensively used in practical applications.

\subsubsection{Computing a Subgradient}
\label{sec:computing_subgradient}
A dual subgradient can be obtained by first minimizing $L(\cdot,\lambda_k)$ and then evaluating an $x _k \in \arg \min_{x \in X} L(x,\lambda_k)$ on the constraints, \emph{i.e.}, $\partial h(\lambda_k)  = g(x_k)$. 
Note that minimizing $L(\cdot,\lambda_k)$ for a fixed $\lambda_k \in \R^m_+$ is an unconstrained convex optimization that can be carried out with a variety of methods, and the choice of using one method or another will depend on the assumptions made on the objective function and constraints. Sometimes it is not possible to exactly minimize the Lagrangian and an approximation is obtained instead, \emph{i.e.}, an $x_k \in X$ such that $L(x_k,\lambda_k) - h(\lambda_k) \le \xi$ where $\xi\ge 0$. This can be equivalently regarded as exactly minimizing the Lagrangian when an approximate Lagrange multiplier is used instead of the true Lagrange multiplier (see Appendix \ref{appendix:preliminaries} for a detailed explanation). That is, we obtain an $x_k \in \arg \min_{x \in X} L(x,\mu_k)$ where $
\mu_k = \lambda_k + \epsilon_k$ and $\epsilon_k \in \R^m$ such that $\mu_k \succeq 0$ ($\epsilon_k$ can be regarded as a perturbation or error in the Lagrange multiplier). 
%

\subsubsection{Convergence}
A standard assumption made to prove the convergence of the subgradient method is that the subdifferential of $h$ is bounded for all $\lambda \succeq 0$. We can ensure this by making the following assumption.
\begin{assumption}
\label{th:boundedset}
$X$ is bounded.
\end{assumption}
\noindent
Observe that since we always have that $\partial h(\lambda)  = g(x)$ for some $x \in X$, if $g(x)$ is bounded for every $x \in X$ then the subgradients of the Lagrange dual function are also bounded. That is, we have that
\begin{align*}
\| \partial h(\lambda)\|_2 \le  \max_{x \in X} \| g(x)  \|_2 : = \sigma_g ,
\end{align*}
and $\sigma_g$ is finite because $g$ is a closed convex function (and so continuous) and $X$ is bounded.

The basic idea behind the convergence of the dual subgradient method with constant step size is that
\begin{enumerate}
\item[(i)] the Euclidean distance between $\lambda_k$ and a vector $\lambda^\star \in \Lambda^\star := \arg \max_{\lambda \succeq 0} h(\lambda)$ decreases \emph{monotonically} when $\lambda_k$ is sufficiently ``far away'' from $\Lambda^\star$;\footnote{Under the Slater condition $\Lambda^\star$ is a bounded subset from $\R^m_+$ (Lemma 1 in \cite{NO09}).}
\item[(ii)] when $\lambda_k$ is sufficiently close to $\Lambda^\star$, it remains in a ball around it. 
\end{enumerate}
Important characteristics of the dual subgradient method are that the size of the ball to which $\lambda_k$ converges depends on $\alpha$; 
 that $\lambda_k$ converges to an $\alpha$-ball around $\Lambda^\star$ in finite time; and that by selecting $\alpha$ sufficiently small we can make the $\alpha$-ball arbitrarily small. 
It is important to highlight as well that the monotonic convergence of $\lambda_k$ to a ball around $\Lambda^\star$ does not imply that the value of the Lagrange dual function improves in each iteration.\footnote{By monotonic convergence we mean the Euclidean distance between $\lambda_k$ and a point in $\Lambda^\star$ decreases.} Yet, since from Assumption \ref{th:boundedset} we have that the Lagrange dual function is Lipschitz continuous\footnote{For any $\lambda_1,\lambda_2 \succeq 0$ we have that $| h(\lambda_1) - h(\lambda_2)| \le \| \lambda_1 - \lambda_2\|_2 \sigma_g$.} for $\lambda \succeq 0$, if the RHS in
\begin{align*}
| h(\lambda_k) - h(\lambda^\star) |  \le \| \lambda_k - \lambda^\star \|_2  \sigma_g
\end{align*}
decreases, then $h(\lambda_k)$ will eventually approach $h(\lambda^\star)$.
%

\section{Subgradient Method with Perturbations}
\label{sec:main}

In this section, we introduce the framework that will allow us to tackle optimization problems with discrete control actions. We begin by considering the following convex optimization problem $\P(\delta)$
\begin{align*}
\begin{array}{ll}
\underset{x \in X}{\text{minimize}} & f(x)\\
\mbox{\text{subject to}} &  g(x) + \delta \preceq 0 
\end{array}
\end{align*}
where $\delta \in \R^m$.
If perturbation $\delta$ were known we could use an interior point method or similar to solve the problem.  However, we will assume that $\delta$ is \emph{not known} in advance and tackle the problem using a Lagrange relaxation on the constraints. 
The interpretation of perturbation $\delta$ will depend on the details of the problem being considered, for example, in a packet switched network $\delta$ may be the (unknown) mean packet arrival rate. 

Because the solution of problem $\P(\delta)$ depends on $\delta$, it will be convenient to define $X_0(\delta) := \{ x \in X  \mid g(x) + \delta \preceq 0\}$, $f^\star (\delta) := \min_{x \in X_0(\delta)} f(x)$, and $x^\star (\delta)$ to be a solution of problem $\P(\delta)$. 
Similarly, we also parameterize the Lagrangian $L(x,\lambda,\delta) = f(x) + \lambda^T  (g(x) + \delta)$, the Lagrange dual function $h(\lambda,\delta) : = \inf_{x \in X} L(x,\lambda,\delta)$, and define the Lagrange dual problem $\D(\delta)$ 
\begin{align*}
\begin{array}{ll}
\underset{\lambda \succeq 0}{\text{maximize}} & h(\lambda,\delta) 
\end{array}
\end{align*}
where $\lambda^\star(\delta)$ is a vector in the set of dual optima $\Lambda^\star (\delta) :=  \arg \max_{\lambda \succeq 0} h(\lambda,\delta)$. 

\subsection{Subgradient Method with Perturbations}
The general version of the subgradient method we consider is the following
\begin{align}
\lambda_{k+1} & =  [\lambda_{k} + \alpha \partial h(\mu_k, \delta_k) ]^+ \label{eq:lp_update} \\
& =  [\lambda_{k} + \alpha (g(x_k) + \delta_k)  ]^+ \notag 
\end{align}
for $k=1,2,\dots$ with $\lambda_1 \in \R^m_+$ and where 
$$x_k \in \arg \min_{x \in X} L(x,\mu_k,\delta_k),$$ $\delta_k \in \R^m$ and $\mu_k = \lambda_k + \epsilon_k$ with $\epsilon_k \in \R^m$.  We will refer to $\delta_k$ and $\epsilon_k$ as perturbations.  
Since parameter $\delta$ is not known in the optimization we have replaced it with a surrogate $\delta_k$, which can be regarded as an approximation or perturbed version of parameter $\delta$. Later, we will add assumptions on the properties that $\delta_k$ must have in order that update (\ref{eq:lp_update}) solves problem $\P(\delta)$.  
An important observation is that for any $\mu_k \in \R^m_+$
\begin{align*}
\arg \min_{x \in X} L(x, \mu_k,\delta)  
& = \arg \min_{x \in X} \{ f(x) + \mu_k^T  (g(x) + \delta)\} \\
& = \arg \min_{x \in X} \{ f(x) + \mu_k^T  g(x)\} 
\end{align*}
and so $g(x_k)$ or the ``partial'' subgradient of $h(\lambda_k,\delta)$ can be obtained independently of perturbation $\delta$.

We are now in the position to present the following lemma. 

\begin{lemma}[Dual Subgradient Method]
\label{th:perturbed_subgradient_method}
Consider the optimization problem $\P(\delta)$ and update (\ref{eq:lp_update}) with $\mu_k \succeq 0$ for all $k$. 
Suppose $\{ \delta_k \}$ is a sequence of points from $\R^m$ such that $\lim_{k \to \infty} k^{-1}\sum_{i=1}^k \delta_i = \delta$.
Then,
\begin{align}
& - \frac{ \| \lambda_1 - \theta\|_2^2 }{2 \alpha k} - \Gamma
 \le   \frac{1}{k} \sum_{i=1}^k  h( \lambda_i,\delta)  - h(\theta,\delta) \label{eq:subgradient_bound}
\end{align}
where $\theta$ is any vector from $\R^m_+$, $\Gamma:=\alpha(\Gamma_a+\Gamma_b+\Gamma_c)+\Gamma_d+\Gamma_e$ and
\begin{align*}
\Gamma_a&:=\frac{1}{2k} \sum_{i=1}^k \|  g(x_i) + \delta \|_2^2,\quad
 \Gamma_b:=\frac{1}{2k} \sum_{i=1}^k\| \delta_i - \delta \|_2^2\\
 \Gamma_c&:= \frac{1}{k}\sum_{i=1}^k (\delta_i - \delta)^T  (g(x_i) + \delta)\\
  \Gamma_d&:=\frac{2}{k} \sum_{i=1}^k\| \epsilon_i \|_2 \| g(x_i) +\delta\|_2\\
 \Gamma_e&:= \frac{1}{k}\sum_{i=1}^k (\lambda_i-\theta)^T  (\delta_i- \delta)
\end{align*}
\end{lemma}
\begin{IEEEproof}
See Appendix \ref{appendix:a}. 
\end{IEEEproof}

Lemma \ref{th:perturbed_subgradient_method} is expressed in a general form and establishes a lower bound on $k^{-1} \sum_{i=1}^k h(\lambda_i,\delta) - h(\theta,\delta)$, where $\theta$ is any vector from $\R^m_+$. When $\theta = \lambda^\star(\delta)$ we can upper bound (\ref{eq:subgradient_bound}) by zero since $h(\lambda,\delta) \le h(\lambda^\star(\delta),\delta)$ for all $\lambda \in \R^m_+$. 
The bound on the left-hand side of (\ref{eq:subgradient_bound}) depends on the properties of perturbations $\delta_k$ and $\epsilon_k$. Next, we provide a detailed analysis of the lower bound and how the different assumptions on the properties of the perturbations result in different types of convergence. A summary of the results can be found in Section \ref{sec:summary}.

\subsubsection{Analysis of Perturbations}
Firstly, by Assumption \ref{th:slater_general} and Lemma 1 in \cite{NO09} we have that $\lambda^\star( \delta)$ is a bounded vector from $\R^m_+$ and therefore the first term on the LHS of (\ref{eq:subgradient_bound})
\begin{align}
\frac{ \| \lambda_1 - \lambda^\star(\delta)\|_2^2 }{2 \alpha k} \label{eq:first_term_alphak}
\end{align} 
is bounded and goes to zero as $k \to \infty$. 
Since (\ref{eq:first_term_alphak}) is divided by $\alpha$,  the convergence rate is inversely proportional to the step size used.

Turning now to $\alpha(\Gamma_a+\Gamma_b+\Gamma_c)$, when $\Gamma_a$, $\Gamma_b$ and $\Gamma_c$ are bounded above then this expression can be made arbitrarily small by selecting $\alpha$ sufficiently small. 
Term $\Gamma_a$ is the sum of the dual subgradients.  If we assume that $X$ is bounded (\emph{c.f.}, Assumption \ref{th:boundedset}) we have that $\| g(x_i) + \delta \|_2$ is upper bounded by some constant $\sigma_g$,\footnote{Note that $g(x_i) + \delta$ is a subgradient of $h(\lambda_k,\delta)$.} and so $\Gamma_a$ is bounded by $\sigma_g^2/2$ for all $i = 1,\dots,k$.
The bounds on terms $\Gamma_b$ and $\Gamma_c$ depend on the characteristics of sequence $\{ \delta_k\}$. 
We consider two cases. Case (i) $\delta_k$, $k=1,2,\dots$ are uniformly bounded. Then, $\Gamma_b$ is trivially uniformly upper bounded for all $k$; and since $(\delta_k - \delta)^T  (g(x_k) + \delta) \le \| \delta_k - \delta \|_2  \| g(x_k) + \delta \|_2 \le \| \delta_k - \delta \|_2 \sigma_g$ by Cauchy-Schwarz, we have that $\Gamma_c$ is also uniformly upper bounded.
Case (ii) $\delta_k$, $k=1,2,\dots$ is a realization of independent random variables with finite variance and kurtosis, but they do not necessarily have to have bounded support. In this case, we can upper bound $\Gamma_b$ and $\Gamma_c$ with probability one asymptotically as $k \to \infty$ using Hoeffding's inequality \cite{Hoe63}. 
Hoeffding's bound can be applied to $\Gamma_b$ directly, and for $\Gamma_c$ it is sufficient to note that 
$  - k^{-1} \sum_{i=1}^k  \sum_{j=1}^m (\delta_i(j) - \delta(j)) (g_l(x_k) + \delta(j))
  \ge - \sum_{j=1}^m | k^{-1}  \sum_{i=1}^k (\delta_i(j) - \delta(j)) | \sigma_g
$
where $\delta(j)$ is the $j$'th component of vector $\delta \in \R^m$, and $l \in \{1,\dots,m\}$. 

Finally, we have the terms $\Gamma_d$ and $\Gamma_e$ which are not scaled by $\alpha$ in the bound.
Since $\Gamma_d$ depends on sequence $\epsilon_k$, the boundedness of the term will depend on the assumptions we make on this perturbation. We consider three cases. 
Case (i) 
$ \| \epsilon_k \|_2 \le \epsilon$ for all $k$ for some $\epsilon > 0$. In this case we have that $\Gamma_d$ can be uniformly upper bounded by $2 \epsilon \sigma_g$. 
Case (ii) $\lim_{k\to \infty} k^{-1} \sum_{i=1}^k  \| \epsilon_i \|_2 =  \epsilon$. We cannot say anything about $\Gamma_d$ for finite $k$, but we will have that $\Gamma_d$ is upper bounded by $2 \epsilon \sigma_g$ when $k \to \infty$. 
An interesting observation is that if $\{ \epsilon_k \}$ were a stochastic process, it would not need to have finite variance in order for $\lim_{k \to \infty} k^{-1} \sum_{i=1}^k \| \epsilon_i \|_2$ to exist and be finite (note this is in marked contrast to perturbation $\delta_k$, which always has to have finite variance).
Case (iii) $\{\|\epsilon_k \|_2 \}$ is a realization of independent random variables with finite variance and mean $\epsilon$. In this case we can use Hoeffding's inequality to give a bound on $\Gamma_d$ with probability one asymptotically as $k \to \infty$.

Term $\Gamma_e$ is perhaps the term for which the analysis is more delicate. In the deterministic subgradient method we have that $\delta_k = \delta$ for all $k$ and so the term is equal to zero for all $k$. 
Observe that when $\Gamma_e$
is nonnegative, then we can ignore the term since this would still leave a lower bound on the LHS of (\ref{eq:subgradient_bound}). 
However, since $\lambda^\star(\delta)$ is not known (we only know it is finite), it is not possible to determine the sign of $\Gamma_e$, and so the term could be unbounded below when $k \to \infty$. 
Because of all this, we will usually require that $\{ \delta_k \}$ is an ergodic process (\emph{i.e.}, $\E(\delta_k) = \delta$ for all $k$) and make use of the fact that $\lambda_k$ and $\delta_k$ are independent for all $k$; in which case 
$
 \E (k^{-1} \sum_{i=1}^k (\lambda_i-\lambda^\star( \delta ))^T  (\delta_i- \delta) ) 
=  k^{-1} \sum_{i=1}^k \E (\lambda_i-\lambda^\star( \delta ))^T  \E(\delta_i- \delta)  
= 0
$
and the \emph{expected} value of the lower bound in Lemma \ref{th:perturbed_subgradient_method} does not depend on term $\Gamma_e$. 

\subsubsection{Summary of the Different Types of Convergence}
\label{sec:summary}
Table \ref{table:summary} provides a summary of the convergence properties of the lower bound obtained in Lemma \ref{th:perturbed_subgradient_method} under the assumption that $\{ \delta_k \}$ is an ergodic stochastic process\footnote{Hence, $\E(\delta_k) = \delta$ for all $k$ where $\delta$ is the perturbation on the constraints given in problem $\P(\delta)$.}  and expectation is taken with respect to random variable $\delta_k$ (since we need term $\Gamma_e$ to vanish). 
By deterministic convergence we mean that it is possible to obtain a lower bound of $ k^{-1} \sum_{i=1}^k  h( \lambda_i,\delta)  - h(\lambda^\star(\delta),\delta) $ for every $k=1,2,\dots$; by w.h.p. that a lower bound can be given with high probability for $k$ large enough; and by $k\to \infty$ that the lower bound will only hold asymptotically, \emph{i.e.}, it is not possible to say anything about the bound for finite $k$. 
\begin{table}
\renewcommand{\arraystretch}{1.5}
\caption{Summary of the convergence of the bound in Lemma \ref{th:perturbed_subgradient_method} depending on the assumptions made on the perturbations ($\E(\delta_k) = \delta$ for all $k$, $\sigma_\delta^2$ is the variance of $\delta_k$, and $\sigma_\epsilon^2$ the variance of $\epsilon_k$).  }
\begin{center}
\begin{tabular}{|l|c|c|c|}
    \cline{2-4}
    \multicolumn{1}{c|}{} & $\sigma_\delta^2=0$ & $\sigma_\delta^2 < \infty$ &  $\sigma_\delta^2 = \infty$ \\ \hline   
$ \|  \epsilon_k \|_2 < \epsilon $ & deterministic & w.h.p. & -  \\ \hline
$\sigma^2_\epsilon < \infty$ & w.h.p.  & w.h.p. & -  \\ \hline
$\sigma^2_\epsilon = \infty$ & $k \to \infty$  & $k \to \infty$ & -  \\ \hline 
\end{tabular}
\end{center}
\label{table:summary}
\end{table}

\subsection{Recovery of Primal Solutions}

We are now in the position to present one of our main theorems, which establishes the convergence of the objective function to a ball around the optimum, and provides bounds on the amount of constraint violation.

\begin{theorem}[Convergence]
\label{th:perturbed_theorem}
Consider problem $\P(\delta)$ and update
\begin{align}
\lambda_{k+1} & =  [\lambda_{k} + \alpha \partial h(\mu_k, \delta_k) ]^+ \label{eq:th_l_udpate}
\end{align}
where $\mu_{k} = \lambda_k + \epsilon_k$ with  $\lambda_1 \in \R^m_+$ and $\{ \epsilon_k \}$ a sequence of points from $\R^m$ such that $\mu_k \succeq 0$ for all $k$.
Suppose $X_0(\delta)$ has nonempty relative interior (the Slater condition is satisfied) and that $ \delta_k $ is an ergodic stochastic process with expected value $\delta$ and $\E(\| \delta_k - \delta \|_2^2) = \sigma_\delta^2$ for some finite $\sigma_\delta^2$.
Further, suppose that $\lim_{k \to \infty} k^{-1} \sum_{i=1}^k \| \epsilon_i \|_2 = \epsilon$ for some $\epsilon \ge 0 $ and that Assumption \ref{th:boundedset} holds. 
Then, 
\begin{flalign*}
 \textup{(i)}  & \
 \E \left( f(\bar x_k) - f^\star(\delta) \right) \le  \frac{\alpha \Theta}{2} + \frac{\| \lambda_{1} \|_2^2}{2\alpha k} +  \frac{2}{k} \sum_{i=1}^k \| \epsilon_i \|_2 \sigma_g  \\
 \textup{(ii)}  & \   \E \left( f(\bar x_k) - f^\star(\delta) \right)
 \ge - 
\frac{\Omega +  \| \E(\bar  \lambda_k) \|_2  (  \| \lambda^\star(\delta) \|_2 + \sqrt{\Omega}     )}{\alpha k}\\
 \textup{(iii)}  
& \  \left\| \E \left( g(\bar x_k) + \delta \right) \right\|_2 \le \frac{1}{\alpha k} \left( \| \lambda^\star(\delta) \|_2 +  \sqrt \Omega \right) \\
  \textup{(iv)}  &  \ \left\| \E \left( \bar \lambda_k \right) \right\|_2 \le \frac{1}{\upsilon}\left(  f(\hat x) - h(\lambda^\star(\delta)) + \frac{\Omega}{\alpha k} \right)
\end{flalign*}
where $\bar x_k = k^{-1} \sum_{i=1}^k x_i$, $\bar \lambda_k : =  k^{-1} \sum_{i=1}^k \lambda_i$, $\Theta : = \sigma^2_g + \sigma_\delta^2$, $\hat x$ is a Slater point (\emph{i.e.}, $g(\hat x) \prec 0$), $\upsilon : = \min_{j \in  \{1,\dots,m\}} -g_j(\hat x)$, and 
$$
\Omega = \| \lambda_1 -  \lambda^\star(\delta)\|_2^2  + \alpha^2  \Theta k  +  2 \alpha \sum_{i=1}^k \| \epsilon_i \|_2 \sigma_g.
$$
\end{theorem}

\begin{IEEEproof}
See Appendix \ref{appendix:a}.
\end{IEEEproof}
Claims (i) and (ii) in Theorem \ref{th:perturbed_theorem} establish that $\E(f(\bar x_k))$ converges to a ball around $f^\star(\delta)$, where the type of convergence will depend on the assumptions made on the perturbations, as indicated in Table \ref{table:summary}. Also, note that choosing $\lambda_1 = 0$ is always a good choice to obtain a sharper upper bound. 
Claim (iii) provides a bound on the expected value of the constraint violation, and claim (iv) says that the expected value of the running average of the Lagrange multipliers is bounded for all $k$. As we will show in Section \ref{sec:remarks}, this will play an important role in establishing the stability of a queueing system. Finally, note that the bounds in claims (ii)-(iv) depend on $\lambda^\star(\delta)$, which is usually not known in the optimization. Fortunately, by Lemma 1 in \cite{NO09}, one can obtain an upper bound on $\lambda^\star (\delta)$, which is enough. Also, observe that we can use the fact that $-h(\lambda^\star(\delta)) \le - h(\lambda)$ for all $\lambda \succeq 0$ and obtain a looser bound  in claim (iv).

Theorem \ref{th:perturbed_theorem} establishes the convergence properties of the dual subgradient update (\ref{eq:th_l_udpate}) under $\delta_k$ and $\epsilon_k$ perturbations without connecting with any specific application. However, by appropriate definition of these perturbations a wide range of situations can be encompassed.  
For example, we can use perturbations $\delta_k$ to relax the perfect knowledge of the constraints, and perturbations $\epsilon_k$ to capture asynchronism in the primal updates (see \cite{VL15}). However, of particular interest here is that we can use the $\epsilon_k$ perturbations in this framework to analyze the use of discrete-valued control actions for solving optimization problem $\P(\delta)$.   We consider this in detail in the next section.


\section{Discrete Control Actions}
\label{sec:actions}

In this section, we present the second main contribution of the paper: how to use $\epsilon_k$ perturbations to equip the dual subgradient method with discrete control actions.
Discrete actions or decisions are crucial in control because many systems are restricted to a finite number of states or choices.
For instance, a traffic controller has to decide whether a traffic light should be red or green, or a video application which streaming quality to use \emph{e.g.}, $\{ 360p, 480p, 720p, 1080p \}$. 

We present two classes of discrete control policies. One that selects discrete control actions in an online or myopic manner based only on the current state of the system, and another batch approach that chooses discrete control actions in blocks or groups.
The latter class is particularly useful for problems where there are constraints or penalties associated with selecting subsets of actions.
For example, in video streaming where the application wants to maximize the quality of the video delivered but at the same time minimize the variability of the quality, and so has constraints on how often it can change the quality of the video stream. 

\subsection{Problem Formulation}
\label{sec:seq_actions}

\begin{figure}
\centering
\begin{tabular}{cc}
\includegraphics[width=0.45\columnwidth]{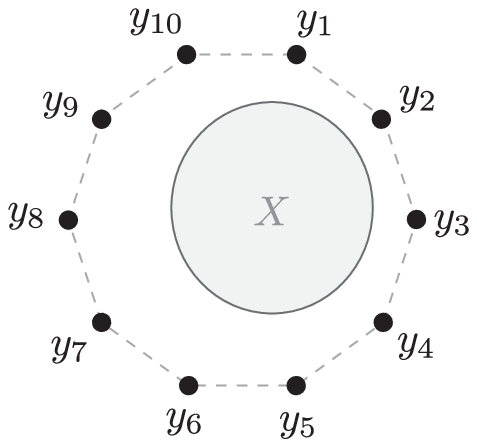} & 
\includegraphics[width=0.45\columnwidth]{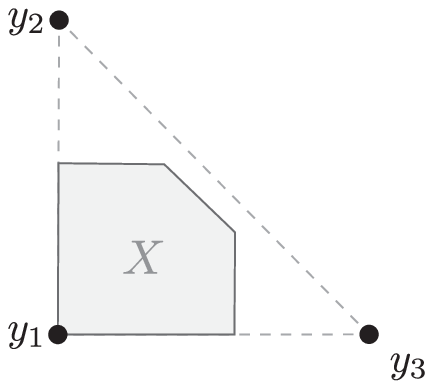}\\
(a) & (b)
\end{tabular}
\caption{Illustrating two action sets $Y$ consisting of a finite collection of points from $\R^2$, and two convex sets $X \subseteq \conv{(Y)}$. The convex hulls of $Y$ are marked in dashed lines.}
\label{fig:setxyconvi}
\end{figure}

We start by introducing the following definitions:
\begin{definition}[Finite Action Set]
$Y$ is a \emph{finite} collection of points from $\R^n$. 
\end{definition}
\begin{definition}[Convex Action Set] $X \subseteq \conv{(Y)}$ and convex.
\end{definition}
We will regard selecting a point $y$ from finite set $Y$ as taking a discrete control action.  The physical action associated with each point in set $Y$ will depend  on the context in which the optimization is applied, \emph{e.g.}, the points in $Y$ may correspond to the actions of setting a traffic light to be red and green.
Figure \ref{fig:setxyconvi} shows an example of two sets $Y$ and respective convex sets $X \subseteq \conv{(Y)}$.

We consider problem $\P(\delta)$ from Section \ref{sec:main}, but now require that the inequality constraints are linear, \emph{i.e.}, $g(x) := Ax$ with $A \in \R^{m \times n}$. The reason for this is the following lemma, which is a restatement of \cite[Proposition 3.1.2]{Mey08}.

\begin{lemma}[Queue Continuity] 
\label{th:skorokhodcontinuity}
Consider updates 
\begin{align*}
\lambda_{k+1} & = [\iter{\lambda}{k} +  \alpha( A x_k + \delta_k)]^+  \\
\mu_{k+1}  & = [ \mu_k +  \alpha (A  y_k + \delta_k)]^+ 
\end{align*}
where $\lambda_1  = \mu_1\ge 0$, $\alpha > 0$, $A \in \R^{m \times n}$, $\delta_k \in \R^m$,  and $\{x_k\}$ and $\{ y_k\}$ are, respectively, sequences of points from $X$ and $Y$.  Suppose $\| \sum_{i=1}^k x_i - y_i \|_2 \le \psi$ for all $k$ and some $\psi \ge 0$.
Then,
\begin{align*}
\| \lambda_k - \mu_k \|_2 \le \epsilon:=2 \alpha \| A \|_2 \psi \qquad k=1,2,\dots
\end{align*}
\end{lemma}

Lemma \ref{th:skorokhodcontinuity} says that when the difference $\| \sum_{i=1}^k x_i - y_i \|_2 $ is uniformly bounded by some constant $\psi$ then $\mu_k$ is an approximate Lagrange multiplier.  And by selecting $x_k$ according to (\ref{eq:lp_update})  we can immediately apply Theorem \ref{th:perturbed_theorem} to conclude that sequence $\{ y_k\}$ of discrete actions (approximately\footnote{Note that if $\psi$ is finite we can make $\epsilon$ (and so $\E(f(\bar x_k) - f^\star(\delta))$ in Theorem \ref{th:perturbed_theorem}) arbitrarily small by selecting step size $\alpha$ sufficiently small. }) solves problem $\P(\delta)$.  This is a key observation.  Not only does it (i) establish that we can solve $\P(\delta)$ using only discrete actions and (ii) give us a testable condition, $\| \sum_{i=1}^k x_i - y_i \|_2\le \psi$, that the discrete actions must satisfy, but it also (iii) tells us that \emph{any} sequence $\{ y_k\}$ of discrete actions satisfying this condition solves $\P(\delta)$.  We therefore have considerable flexibility in how we select actions $y_k$.  In other words, we have the freedom to select from a range of different optimal control policies without changing the underlying convex updates.  One way to use this freedom is to select an optimal control policy that satisfies specified constraints, \emph{e.g.}, that does not switch traffic lights between green and red too frequently or which minimizes the use of ``costly'' actions.   Such constraints are often practically important yet are difficult to accommodate within classical control design frameworks.

With this in mind, in the rest of this section we consider methods for constructing sequences $\{ y_k\}$ of discrete actions that stay close to a sequence $\{ x_k\}$ of continuous-valued updates in the sense that $\| \sum_{i=1}^k x_i - y_i \|_2$ is uniformly bounded.  We begin by establishing that for any sequence $\{x_k\in X \subseteq \conv(Y)\}$ there always exists discrete-valued sequences $\{y_k\in Y\}$ such that $\| \sum_{i=1}^k x_i - y_i \|_2$ is uniformly bounded.

\subsection{Existence of Discrete Sequences}
\label{sec:existence}
It will prove convenient to exploit the fact that each point $x \in X$ can be written as a convex combination of points from $Y$. Collect the points in $Y$ as columns in matrix $W$ and define
\begin{align*}
E : = &  \{v_1,\dots,v_{|Y|}\}, \\
U := & \conv{(E)} = \{ u \in [0,1]^{|Y|} \mid \1^T  u  =1 \},
\end{align*}
where $v_j$ is an $|Y|$-dimensional standard basis vector, \emph{i.e.},  all elements of vector $v_j$ are equal to $0$ except the $j$'th element that is equal to $1$, and $U$ is the $|Y|$-dimensional simplex.  Since we can always write a vector $x_i \in X$ as the convex combination of points from $Y$ there exists at least one vector $u_i \in U$ such that $x_i = Wu_i$.\footnote{A vector $u_i$ can be obtained, for example, by solving the optimization problem $\min_{u \in U} \| x_i - W u \|_2^2$. The non-uniqueness of the solution comes from Carath\'eodory's theorem---see, for example, \cite{Ber60}.} Similarly, there exists a vector $e_i \in E$ such that  $y_i = We_i$.
Hence, 
\begin{align*}
\left\| \sum_{i=1}^k x_i - y_i \right\|_2 \!
= \left\| \sum_{i=1}^k W (u_i - e_i) \right\|_2  
\! \le \| W \|_2 \left\| \sum_{i=1}^k u_i - e_i \right\|_2 
\end{align*} 
and therefore showing that $ \| \sum_{i=1}^k u_i - e_i \|_2$ is uniformly bounded is sufficient to establish the boundedness of $ \| \sum_{i=1}^k x_i - y_i \|_2$.

We have the following useful lemma.
\begin{lemma}
\label{th:blocksequence}
Let $E$ be a set containing the $|Y|$-dimensional standard basis vectors,  
$U : = \conv{(E)}$, and 
$D := \{ \delta \in \R^{|Y|} \mid \delta^T  \1 = 0, \| \delta \|_\infty \le 1 \}$. 
For any vector $\delta \in D$, and sequence $\{ u_i \}_{i=1}^{|Y|}$ of points from $U$, there exists at least one sequence $\{ e_i\}_{i=1}^{|Y|}$ of points from $E$ such that 
\begin{align}
\left( \delta + z - z'  \right) \in D,
\end{align}
where $z := \sum_{i=1}^{|Y|} u_i$ and $z' := \sum_{i=1}^{|Y|} e_i$.   That is, $\1^T  ( \delta + z - z') = 0$ and $\|  \delta + z - z' \|_\infty \le 1$. 
\end{lemma}
\begin{IEEEproof}
See Appendix \ref{appendix:c}.
\end{IEEEproof}
Let $\{ u^{(1)}_{i} \}_{i=1}^{|Y|}$ be a sequence of points from $U$, and $\{ e^{(1)}_i \}_{i=1}^{|Y|}$ a sequence of points from $E$ such that $(z_1 - z'_1) \in D$ where $z_1=\sum_{i=1}^{|Y|} u^{(1)}_i$ and $z'_1= \sum_{i=1}^{|Y|} e^{(1)}_i$.  By Lemma \ref{th:blocksequence} such a sequence $\{ e^{(1)}_i \}_{i=1}^{|Y|}$ always exists.
Similarly, for another sequence $\{ u^{(2)}_i \}_{i=1}^{|Y|}$ of points from $U$, we can construct a sequence $\{ e^{(2)}_i \}_{i=1}^{|Y|}$ of points from $E$ such that $ (z_2 - z_2' + (z_1 - z'_1)) \in D$ where $z_2$ and $z'_2$ are, respectively, the sum of the elements in sequences $\{ u^{(2)}_i \}_{i=1}^{|Y|}$ and $\{ e^{(2)}_i \}_{i=1}^{|Y|}$.
Repeating, it follows that for sequences $\{ u^{(\tau)}\}_{i=1}^{|Y|}$, $\tau \in \{1,2,\dots, K\}$ we can construct sequences $\{ e^{(\tau)}\}_{i=1}^{|Y|}$, $\tau \in \{1,2,\dots, K\}$ such that
\begin{align*}
& (((\cdots (( z_1 - z'_1) + z_2 - z'_2 )+ \dots + z_{K-1} - z'_{K-1}) \notag \\
& \quad  + z_K - z_K')  ) =  \left( \sum_{\tau=1}^K z_\tau - z'_\tau \right) \in D,
\end{align*} 
where $z_\tau$ and $z'_\tau$ are the sum of the elements in the respective sequences.  It follows that $\| \sum_{i=1}^k u_i -e_i \|_2 \le \sqrt{|Y|}$ for $k \in \tau |Y|$, $\tau \in \Z_+$ and $\| \sum_{i=1}^k u_i -e_i \|_2 \le \sqrt{|Y|(1+2|Y|)^2}$ for all $k=1,2,\dots$, since the sequences can diverge element-wise by at most $2|Y|$ over the $|Y|$ steps between $\tau |Y|$ and $(\tau+1) |Y|$.  We therefore have the following result.
\begin{theorem}[Existence of Discrete Sequences]\label{cor:existence}
For any sequence $\{ u_k\}$ of points from $U$ there exists a sequence $\{ e_k\}$ of points from $E$ such that $\| \sum_{i=1}^k u_i -e_i \|_2$ is uniformly bounded for all $k \in \{1,2,\dots\}$.
\end{theorem}

Note that since we can always permute the entries in sequence $\{ e_k\}$ while keeping $\| \sum_{i=1}^k u_i -e_i \|_\infty$ bounded, the existence of one sequence implies the existence of many (indeed exponentially many since the number of permutations grows exponentially with the permutation length).

\subsection{Constructing Sequences of Discrete Actions Using Blocks}
\label{sec:seq_block_actions}
We now present our first method for constructing sequences of discrete actions, which uses a block-based approach.  
\begin{lemma}
\label{th:magic_update_lemma}
Consider the setup of Lemma \ref{th:blocksequence} and select
\begin{align*}
e_i \in \arg \min_{e \in E} 
 \| ( \delta + z - \sum_{\kappa=1}^{i-1} e_\kappa ) - e \|_\infty
 ,\ i = 1,\dots,|Y|
\end{align*}
where $z := \sum_{i=1}^{|Y|} u_i$ and $\delta\in D$.  Then, $- \1 \preceq \delta + z - z' \preceq \1$, with $z' = \sum_{i=1}^{|Y|} e_i$. 
\end{lemma}
\begin{IEEEproof}
See Appendix \ref{appendix:c}.
\end{IEEEproof}
Partitioning sequence $\{u_k\}$, $k \in \{1,2,\dots\}$, of points from $U$ into subsequences $\{ u_i^{(\tau)} \}_{i=1}^{|Y|}$, $\tau\in \Z_+$ with $u_i^{(\tau)}=u_{\tau|Y|+i}$ and applying Lemma \ref{th:magic_update_lemma} recursively yields a sequence $\{e_i\}$ such that $\| \sum_{i=1}^k u_i -e_i \|_2$ is uniformly bounded.  That is, we have the following.
\begin{theorem}
Let $\{u_k\}$ be a sequence of points from $U$ partitioned into subsequences $\{ u_i^{(\tau)} \}_{i=1}^{|Y|}$ with $u_i^{(\tau)}=u_{\tau|Y|+i}$, $\tau\in \Z_+$.  For $i \in \{ 1,\dots,|Y|\}$, $\tau \in \{0,1,\dots \}$ select
\begin{align}
e^{(\tau)}_i \in \arg \min_{e \in E} \left\| \left( z_\tau -  \sum_{\kappa=1}^{i-1}  e_\kappa^{(\tau)} \right) - e \right\|_\infty \label{eq:magic_update}
\end{align}
where $z_\tau := \sum_{i=1}^{|Y|} u_i^{(\tau)}$.  Then $\| \sum_{i=1}^k u_i -e_i \|_2$ is uniformly bounded for all $k=1,2,\dots$, where $e_{\tau|Y|+i}=e_i^{(\tau)}$. 
\end{theorem}

Observe that, with this approach, the construction of a subsequence $\{ e_i^{(\tau)} \}_{i=1}^{|Y|}$ requires that subsequence $\{ u_i^{(\tau)} \}_{i=1}^{|Y|}$ is known.  Hence, we refer to this as a block or batch approach.  When sequence $\{u_k\}$ is observed in an online manner, then sequence $\{e_k\}$ must be constructed with a delay of $|Y|$ elements relative to $\{u_k\}$ since in order to construct $\{ e_i^{(\tau)} \}_{i=1}^{|Y|}$ we must wait for $|Y|$ elements until $\{ u_i^{(\tau)} \}_{i=1}^{|Y|}$ is observed.  

Note that by a similar analysis we can immediately generalize this method of construction to situations where we partition sequence $\{u_k\}$, $k \in \{ 1,2,\dots \}$, into subsequences $\{ u_i^{(\tau)} \}_{i=1}^{T_\tau|Y|}$, $T_\tau\in \N$, $\tau\in \Z_+$ with $u_i^{(\tau)}=u_{\sum_{t=0}^\tau T_t|Y|+i}$ \emph{i.e.}, where the subsequences can be different lengths so long as they are all some multiple of $|Y|$.

\subsubsection{Constrained Control Actions}
We can permute sequence $\{ e_i^{(\tau)} \}_{i=1}^{|Y|}$ arbitrarily while keeping $\| \sum_{i=1}^k u_i -e_i \|_2$ uniformly bounded (since the vectors $u_i^{(\tau)}$ and $e_i^{(\tau)}$ both lie in the unit ball then $\|u_i^{(\tau)}-e_j^{(\tau)}\|_\infty\le 2$ for all $i$, $j\in\{1,\dots,|Y|\}$).   When the sequence of admissible actions is constrained, this flexibility can be used to select a sequence of actions which is admissible.   For example, sequence $\{0,1,0,1\}$ might be permuted to $\{0,0,1,1\}$ if the cost of changing the action taken in the previous iteration is high.

\subsection{Constructing Sequences of Discrete Actions Myopically}
\label{sec:online}
We now consider constructing a discrete valued sequence $\{e_k\}$ in a manner which is myopic or ``greedy'', \emph{i.e.}, that selects each $e_k$, $k=1,2,\dots$ by only looking at 
$u_i$, $i=1,\dots,k$ and $e_i$, $i=1,\dots,k-1$.  
We have the following theorem.

\begin{theorem}
\label{th:online_actions}
Let $\{ u_k \}$ be a sequence of points from $U$.  Select
\begin{align}
e_{k} & \in \arg \min_{e\in E} \| s_{k-1} + u_{k} - e \|_\infty, \label{eq:ekupdate}
\end{align}
where $s_k = \sum_{i=1}^k (u_i-e_i)$. 
Then, we have that
$
-\1 \preceq  s_k \preceq (|Y|-1) \1,
$
and 
\begin{align*}
\left\| \sum_{i=1}^k u_i - e_i  \right\|_2 \le C:= \sqrt{|Y|} (|Y| - 1) 
\end{align*}
\end{theorem}
\begin{IEEEproof}
See Appendix \ref{appendix:c}.
\end{IEEEproof}

Theorem \ref{th:online_actions} guarantees that by using update (\ref{eq:ekupdate}) the difference $\| \sum_{i=1}^k u_i - e_i \|_2$ is uniformly bounded. 
Observe that update (\ref{eq:ekupdate}) selects a vector $e \in E$ that decreases the largest component of vector $s_k$, and so it does not provide flexibility to select other actions in $Y$.  That is, the benefit of myopic selection comes at the cost of reduced freedom in the choice of discrete control action.  Regarding complexity, solving (\ref{eq:ekupdate}) in general involves using exhaustive search.  However, it is not necessary to solve (\ref{eq:ekupdate}) for every $u_k$, $k=1,2,\dots$ so long as the difference between the steps when (\ref{eq:ekupdate}) is performed is bounded.  This is because the elements of $u_k$ and $e_k$ can diverge by at most $2$ at every step (recall both vectors lie in the unit ball) and so remain bounded between updates.  Hence, the cost of (\ref{eq:ekupdate}) can be amortized across steps.  We have the following corollary.
\begin{corollary}
\label{th:seq_online}
Consider the setup of Theorem \ref{th:online_actions} and suppose update  (\ref{eq:ekupdate}) is performed at steps $k\in\{\tau_1,\tau_2,\dots \} := \mathcal T \subseteq \N$; otherwise, $e_{k}$ is selected equal to $e_{k-1}$. 
Then, we have that
$
-\bar \tau  \1 \preceq  s_k \preceq \bar \tau  (|Y|-1)  \1 ,
$ for all $k$ where  $\bar \tau = \max_{j \in \{1,2,...\}} \{\tau_{j+1}-\tau_j \}$ 
and 
\begin{align*}
\left\| \sum_{i=1}^k u_i - e_i  \right\|_2 \le \bar \tau C.  
\end{align*}
\end{corollary}

\begin{IEEEproof}
See Appendix \ref{appendix:c}.
\end{IEEEproof}

\subsubsection{Constrained Control Actions}
As already noted, we are interested in constructing sequences of actions in a flexible way that can be easily adapted to the constraints on the admissible actions.   Corollary \ref{th:seq_online} allows us to accommodate one common class of constraints on the actions, namely that once an action has been initiated it must run to completion.  For example, suppose we are scheduling packet transmissions where the packets are of variable size and a discrete action represents transmitting a single bit, then Corollary \ref{th:seq_online} ensures that a sequence of bits can be transmitted until a whole packet has been sent. The condition that $\bar \tau $ must be finite corresponds in this case to requiring that packets have finite length. 

In the case of myopic updates it is difficult to give an algorithm without specifying a problem; nevertheless, we can establish the conditions that a generic algorithm should check when selecting a sequence of actions. As shown in the previous section, for finite $|Y|$ it is possible to construct a sequence of actions that breaks free from the past for a subsequence that is sufficiently large. The same concept can be applied to the myopic case, but now we must ensure that $\| s_k \|_\infty$ is bounded for all $k$. This motivates the following theorem. 
%

\begin{theorem}
\label{th:gen_actions}
Let $\{ u_k \}$ be a sequence of points from $U$. 
For any sequence $\{e_k\}$ of points from $E$ we have that 
\begin{align*}
\left\| \sum_{i=1}^k u_i - e_i  \right\|_2 \le  \mathcal \gamma_k C 
\end{align*}
where $\gamma_k = - \min_{j \in \{1,\dots,|Y| \} } s_k(j)$, $s_k = \sum_{i=1}^k u_i -e_i$ and $C = \sqrt{|Y|} (|Y| - 1)$.
\end{theorem}
\begin{IEEEproof}
See Appendix \ref{appendix:c}.
\end{IEEEproof}

Theorem \ref{th:gen_actions} says that when we can construct a sequence of actions $\{ e_k \}$, such that $\gamma_k$ is bounded then the difference $\| \sum_{i=1}^k u_i - e_i \|_2$ will be bounded.  However, now $e_k$ does not need to be obtained as in (\ref{eq:ekupdate}) as long as it is selected with some ``care''.  Namely, by not selecting actions that decrease lower bound $\gamma_k$ ``excessively''. For example, choosing a vector $e_k$ that decreases a positive component of vector $s_k$ will be enough. 
In addition to providing flexibility to respect constraints in the admissible actions, the implications of this are also important for scalability, namely when action set is  large we do not need to do an exhaustive search over all the elements to select a vector from $E$.

\section{Remarks}
\label{sec:remarks}

\subsection{Discrete Queues}
Perturbations on the Lagrange multipliers can be used to model the dynamics of queues that contain discrete valued quantities such as packets, vehicles or people. Suppose that the vector of approximate Lagrange multipliers $\mu_k$ takes a queue-like form, \emph{i.e.},
\begin{align*}
\mu_{k+1} = [\mu_{k} + \alpha \rho_k]^+ , 
\end{align*}
where  $\rho_k \in Y\subset \Z^m$ and represents the number of discrete elements that enter and leave the queues. Since $[\cdot]^+$ is a homogeneous function we can define $Q_k := \mu_k / \alpha$ and obtain iterate
\begin{align}
Q_{k+1}  = [Q_k + \rho_k]^+. \label{eq:queueudpateex}
\end{align}
That is, since $\rho_k$ takes values from a subset of integers we have that $Q_k$ is also integer valued, and therefore update (\ref{eq:queueudpateex}) models the dynamics of a queue with discrete quantities.  Provided $\epsilon_k=\lambda_k-\alpha Q_k$ satisfies the conditions of Theorem \ref{th:perturbed_theorem} then we can use $\alpha Q_k$ (which is equal to $\mu_k$) in update $(\ref{eq:th_l_udpate})$.

\subsection{Queue Stability}
A central point in max-weight approaches is to show the stability of the system, which is usually established by the use of Foster-Lyapunov arguments. In our approach it is sufficient to use the boundedness of the Lagrange multipliers (\emph{i.e.}, claim (iv) in Theorem \ref{th:perturbed_theorem}) and the fact that the difference $\| \lambda_k - \mu_k \|_2 = \| \lambda_k - \alpha Q_k \|_2$ remains uniformly bounded. 
Observe that the uniform boundedness of $ \| \lambda_k - \alpha Q_k \|_2$ implies that $\E (k^{-1} \sum_{i=1}^k \alpha Q_i )  =  \alpha \E (k^{-1} \sum_{i=1}^k  Q_i )\prec \infty $ for all $k$ as well. 
Also, by the linearity of the expectation we can take the expectation inside the summation, and by dividing by $\alpha$ and taking $k \to \infty$ we can write
\begin{align*}
\lim_{k \to \infty} \frac{1}{k} \sum_{i=1}^k \E(Q_i)  \prec \infty,
\end{align*}
which is the definition of strong stability given in the literature of max-weight---see, for example, \cite[Definition 4]{Nee10}. 

\subsection{Primal-Dual Updates}
The focus of the paper has been on solving the Lagrange dual problem directly, but the analysis encompasses primal-dual approaches as a special case. For instance, instead of obtaining an $x_k \in \arg \min_{x \in X} L(x, \lambda_k,\delta)$ in each iteration we could obtain an $x_k$ that provides ``sufficient descent'' \emph{i.e.}, an update that makes the difference $L(x_k,\lambda_k,\delta) - h(\lambda_k)$ decrease monotonically. 
This strategy is in spirit very close to the dynamical system approaches presented in \cite{Sto05,ES05,VL16}, which usually require the objective function and constraints to have bounded curvature.  
In our case, having that the difference $L(x_k,\lambda_k,\delta) - h(\lambda_k)$ decreases monotonically translates into having a sequence $\{ \epsilon_k\}$ that converges to zero and so the perturbation vanishes.\footnote{See Section \ref{sec:computing_subgradient} for a more detailed explanation of how the difference $L(x_k,\lambda_k,\delta) - h(\lambda_k)$ relates to the $\epsilon_k$ perturbations.} 

\subsection{Unconstrained Optimization}

The main motivation for using Lagrange relaxation is to tackle resource allocation problems of the sort tackled by max-weight approaches in the literature. However, our results (Lemma \ref{th:perturbed_subgradient_method}) naturally extend to unconstrained optimization problems.  In this case $h$ becomes the unconstrained objective function that takes values from $\R^m \to \R$, the Lagrange multiplier $\lambda$ is the ``primal'' variable from a convex set $D$, and $\Pi_D$ (instead of $[\cdot]^+$) the Euclidean projection of a vector $\lambda \in \R^m$ onto $D$. 
The proof of Lemma \ref{th:perturbed_subgradient_method} remains unchanged, and it is sufficient to note that the Euclidean projection onto a convex set is nonexpansive. That is, 
$\| \Pi_\D(\lambda_1) - \Pi_D ( \lambda_2) \|_2 \le \|  \lambda_1 - \lambda_2 \|_2$ for all $\lambda_1, \lambda_2 \in \R^m$. 
%

\section{Numerical Example}
\label{sec:example}

\subsection{Problem Setup}
\begin{figure}
\centering
\includegraphics[width=0.6\columnwidth]{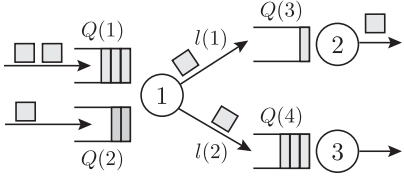}
\caption{Illustrating the network in the example of Section \ref{sec:example}. The access point (node 1) receives packets in queue 1 and queue 2, and sends them, respectively, to nodes 2 and 3. The packets sent by nodes 2 and 3 leave the system.}
\label{fig:ex_net}
\end{figure}

Consider the network shown in Figure  \ref{fig:ex_net} where an Access Point (AP), labelled as node 1, transmits to two wireless stations, labelled as nodes 2 and 3.
Time is slotted and in each slot packets arrive at the queues $Q(1)$ and $Q(2)$ of node 1, which are later transmitted to nodes 2 and 3.  In particular, node 1 transmits packets from $Q(1)$ to $Q(3)$ (node 2) using link $l(1)$, and packets from $Q(2)$ to $Q(4)$ (node 3) using link $l(2)$ (see Figure \ref{fig:ex_net}).
We represent the connection between queues using the incidence matrix $A \in \{-1,0,1 \}^{n \times l}$ where $-1$ indicates that a link is leaving a queue; $1$ that a link is entering a queue; and $0$ that a queue and a link are not connected. For example, a $-1$ in the $j$'th element of the $i$'th row of matrix $A$ indicates that link $j$ is leaving queue $i$. 
The incidence matrix of the network illustrated in Figure \ref{fig:ex_net} is given by 
\begin{align}
A = 
\begin{bmatrix}
-1 & \phantom{-}0 & \phantom{-}1 & \phantom{-}0 \\ 
\phantom{-}0  & -1 & \phantom{-}0 & \phantom{-}1 
\end{bmatrix}^T .
\label{eq:incidence_matrix}
\end{align}

In each time slot, the AP (node 1) takes an action from action set  $Y :=  \{ y{(0)}, y{(1)}, y{(2)}\} = \{ [0,0]^T , [1,0]^T, [0,1]^T \}$, where each action indicates which link to activate. For simplicity of exposition we will assume that activating a link corresponds to transmitting a packet, and therefore selecting action $y(1)$ corresponds to transmitting a packet from $Q(1)$ to node $Q(3)$; action $y(2)$ to transmitting a packet from $Q(2)$ to $Q(4)$; and action $y(0)$ to doing nothing, \emph{i.e.}, not activating any of the links. 

The goal is to design a scheduling policy for the AP (select actions from set $Y$) in order to minimize a convex cost function of the average throughput (\emph{e.g.}, this might represent the energy usage), and ensure that the system is stable, \emph{i.e.}, the queues do not overflow and so all traffic is served.

The convex or fluid-like formulation of the problem is 
\begin{align}
\begin{array}{ll}
 \underset{x \in X}{\text{minimize}} & \quad f(x) \\
 \underset{}{\text{subject to}} & \quad Ax + b \preceq 0
\end{array}
\label{eq:fluid_problem}
\end{align}
where $f: \R^2 \to \R$ is a convex utility function, $X \subseteq \conv{(Y)}$, $A$ the incidence matrix given in (\ref{eq:incidence_matrix}), and $b$ a vector that captures the \emph{mean} exogenous packet arrivals/departures in the system.\footnote{More precisely, $b(1)$ and $b(2)$ capture, respectively, the mean arrival into $Q(1)$ and $Q(2)$; and $b(3)$ and $b(4)$, respectively, the mean departure rate from $Q(3)$ and $Q(4)$.}

\subsection{Unconstrained Control Actions}
Problem (\ref{eq:fluid_problem}) can be solved with the optimization framework presented in Section \ref{sec:main}. That is, with update
\begin{align}
\lambda_{k+1} & = [\lambda_{k} + \alpha (Ax_k + B_k)]^+ \label{eq:subexupdate}
\end{align}
where $x_k  \in \arg \min_{x \in X} \  \{ f(x) + \lambda_k^T A x \}$, and note we have replaced $b$ with random variable $B_k$ in order to capture the fact that packet arrivals at node 1 might be a realization of a stochastic process.

Observe from update (\ref{eq:subexupdate}) that by selecting an $x_k \in  X:=\conv(Y)$ we obtain the fraction of time each of the links should be used in each iteration, but not which packet to transmit from each of the queues in each time slot. Nonetheless, we can easily incorporate (discrete) control actions $y_k\in Y$ by using, for example, Theorem \ref{th:online_actions}.   Also, note that if we define  approximate Lagrange multiplier $\mu_{k+1} = [\mu_k + \alpha(A y_k + B_k)]^+$ and let $Q_k : = \mu_k / \alpha$ we obtain queue dynamics 
\begin{align*}
Q_{k+1} = [Q_k + Ay_k + B_k]^+,
\end{align*}
which are identical to those of the real queues in the system.
By Theorem \ref{th:perturbed_theorem} we can use update $x_k  \in \arg \min_{x \in X} \  \{ f(x) + \mu_k^T A x \}=\arg \min_{x \in X} \  \{ f(x) + \alpha Q_k^T A x \}$, and with this change we now do not need to compute the Lagrange multiplier $\lambda_k$. Looking at the current queue occupancies $Q_k$ in the system is enough for selecting control actions. 

\subsection{Constrained Control Actions} 

We now extend the example to consider the case where the admissible sequence of control actions is constrained. 
Specifically, suppose it is not possible to select action $y{(1)}$ after $y{(2)}$ without first selecting $y{(0)}$; and in the same way, it is not possible to select $y{(2)}$ after $y{(1)}$ without first selecting $y{(0)}$. However, $y{(1)}$ or $y{(2)}$ can be selected consecutively. An example of an admissible sequence would be
\begin{align*}
\{ y{(1)}, \dots, y{(1)}, y{(0)}, y{(2)}, \dots, y{(2)}, y(0),y{(1)},\dots,y(1)\}. 
\end{align*}
This type of constraint appears in a variety of applications ---they are usually known in the literature as reconfiguration or switchover delays \cite{CM15}--- and capture the cost of selecting a ``new'' control action.
In this wireless communication example, the requirement of selecting action $y(0)$ every time the AP wants to change from action $y(1)$ to $y(2)$ and from $y(2)$ to $y(1)$ might be regarded as the time required for measuring Channel State Information (CSI) in order to adjust the transmission parameters.\footnote{The CSI in wireless communications is in practice measured periodically, and not just at the beginning of a transmission, but we assume this for simplicity. The extension is nevertheless straightforward.}

\begin{figure}
\centering
 \includegraphics[width=0.45\columnwidth]{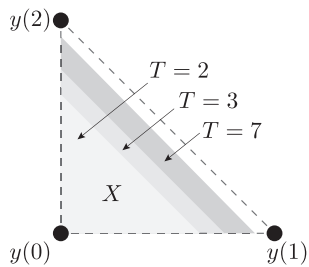} 
\caption{Illustrating how set $X$ defined in (\ref{eq:setX}) changes depending on parameter $T$. The convex hull of $Y$ is indicated with dashed lines. 
}
\label{fig:ex_set}
\end{figure}

In this case, the constraints on the selection of control actions will affect the definition of set $X$ in the problem.  
To see this, observe that if we select a sequence of actions in blocks of length $T|Y|$ with $T \in \N$, and $y{(1)}$ and $y{(2)}$ appear (each) consecutively, then $y{(0)}$ should appear \emph{at least} twice in order for the subsequence to be compliant with the transmission protocol. Conversely, any subsequence of length $T|Y|$ in which $y{(0)}$ appears at least twice can be reordered to obtain a subsequence that is compliant with the transmission protocol. For example, if $T = 3$ and we have subsequence
\begin{align*}
\{ y(1), y(2),  y(1), y(2),  y(2),  y(0),   y(1),  y(2),  y(0)  \},
\end{align*}
we can reorder it to obtain
\begin{align*}
\{ y(0), y(1),  y(1), y(1),  y(0),  y(2),   y(2),  y(2),  y(2)  \},
\end{align*}
which is a subsequence compliant with the transmission protocol. 
Since from Section \ref{sec:seq_block_actions}, we can always choose a subsequence of discrete actions and then reorder its elements, we just need to select a set $X$ such that $y{(0)}$ can be selected twice in a subsequence of $T|Y|$ elements. 
This will be the case when every point $x \in X$ can be written as a convex combination of points from $Y$ that uses action $y({0})$ at least fraction $2/(T|Y|)$. That is, when
\begin{align}
X \subseteq \left(1-  \frac{2}{T|Y|}\right) \, \mathrm{conv}(Y). \label{eq:setX}
\end{align}
Observe from (\ref{eq:setX}) that as $T$ increases we have that $X \to \conv(Y)$, which can be regarded as increasing the capacity of the network since it will be possible to use links 1 and 2 a larger fraction of time. Figure \ref{fig:ex_set} illustrates the capacity of the network (set $X$) changes depending on parameter $T$. 

\subsubsection{Simulation}
We run a simulation using $f(x) = \| Sx \|_2^2$, $S = \mathrm{diag}([1,3])$, $b = [0.25, 0.5,-1,-1]^T$,  $\lambda_1 = \alpha Q_1 = 0$ and $\alpha = 0.01$.
At each iteration we perform update (\ref{eq:subexupdate}) where $x_k \in \arg \min_{x \in X} \{ f(x) + \alpha Q_k^T A x \}$; and $B_k(1)$ and $B_k(2)$ are Bernoulli random variables with mean $b(1)$ and $b(2)$ respectively; $B_k(3)$ and $B_k(4)$ are equal to $-1$ for all $k$ and so the service of nodes 2 and 3 is deterministic. 
Discrete control actions are selected with update (\ref{eq:magic_update})
with $T = 3$, and so we have that the number of elements in a block or subsequence is 9. The elements in a block are reordered in order to match the protocol constraints of the AP. 

\begin{figure}[ht]
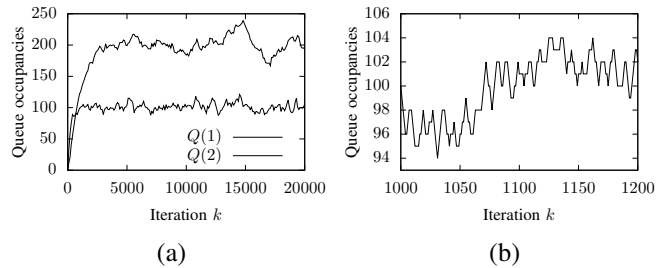

\centering
\setlength{\tabcolsep}{0pt}
\begin{tabular}{cc}
{\resizebox{0.5\columnwidth}{!}{\input{Figures/valls5.tex}}} & 
 {\resizebox{0.5\columnwidth}{!}{\input{Figures/valls6.tex}}} \\
(a) &  (b)
\end{tabular}
\caption{Illustrating the occupancy of queues 1 and 2 in the network. Figure (b) shows the detail of (a) for an interval of 200 time slots.  Observe from (b) that the occupancy of the queues is integer valued. }
\label{fig:queues}
\end{figure}

Figure \ref{fig:queues} plots the occupancy of the queues in the system. Observe that $\alpha Q_k$ converges to a ball around $\lambda^\star = [2,1]^T$, the optimal dual variable associated to the fluid/convex problem (\ref{eq:fluid_problem}). Figure \ref{fig:queues} (b) is the detail of (a) for an interval of 200 iterations, and shows that queues are indeed integer valued.
\begin{figure}[hbt!]
\centering
{\resizebox{0.7\columnwidth}{!}{\input{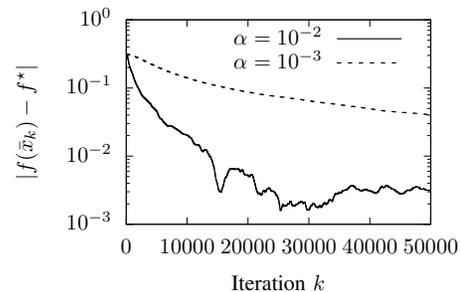}}}
\caption{Illustrating the difference $| f(\bar x_k) - f^\star | $ for different step sizes.}
\label{fig:objective}
\end{figure}

Figure \ref{fig:objective} plots the convergence of the objective function for step sizes $\alpha \in \{ 10^{-2}, 10^{-3} \}$. Observe that using a smaller step size heavily affects the convergence time, which is in line with the accuracy vs. convergence time tradeoff in subgradient methods with constant step size.
%

\section{Related Work}
\label{sec:relatedwork}

In this section, we explain the differences of our contributions with previous work.

\subsection{Contribution \textup{(i)} --- Unifying Framework}

\subsubsection{Lagrange Multipliers and Queues} Concerning the existence of a connection  between the  discrete-valued  queue  occupancy  in  a  queueing  network and  continuous-valued  Lagrange  multipliers, \cite{Sto05} shows that asymptotically, as the design  parameters $\beta\rightarrow 0$ and $t\rightarrow\infty$, the scaled queue occupancy converges to the set of dual optima.  Also, in \cite{Hua11} it is established that a discrete queue update is exponentially attracted to a ``static'' vector of Lagrange multipliers.
Regarding queues with continuous-valued occupancy (so ``fluid'' type queues), it is shown in \cite{LS04} that the scaled queue occupancy is equal to the Lagrange multipliers generated by an associated dual subgradient update. However, this approach does not encompass common situations such as discrete packet arrivals.  In \cite{LS06}, the authors extend this observation to consider queues that contain errors due to asynchronism in the network but do not present any analytic bounds.

Our work extends this by identifying queues with $\alpha$-scaled approximations of the Lagrange multipliers, providing both a non-asymptotic bound on the approximation and establishing how this affects convergence. In particular, we show how the convergence of the dual subgradient method is affected by $\epsilon_k$-subgradients in the form of approximate Lagrange multipliers. 

\subsubsection{Generality}

Our framework is an extension of Nedi\'c's and Ozdaglar's work in \cite{NO09} to consider stochastic and $\epsilon_k$-subgradients, which we later identify, respectively, with myopic and discrete actions. 
Unlike previous works, we have abstracted these two features and made them accessible from within a clean mathematical framework that does not rely on a specific application. Hence, our contribution is in spirit very different to the works in \cite{LS04, LS06, ES06, ES07, CLCD06} which focus on specific congestion control and scheduling applications.  For example, \cite{ES06} addresses the problem of designing  a congestion controller that is more gradual than previous \emph{dual} controllers so that it ``mimics'' TCP's behavior, and \cite{CLCD06} deals with the joint congestion control, scheduling, and routing in networks with time-varying channel conditions.

\subsubsection{Convergence}

The max-weight analysis in the literature \cite{LS04, LS06, ES06, ES07, CLCD06} provides asymptotic convergence and stability results. This is in marked contrast to our framework, in which we can make quantitative statements of the system state in finite time. In particular, we provide upper and lower bounds on the objective function, bound the constraints violation, and show that the expected value of the Lagrange multipliers is bounded. Compare, for example, \cite[Theorem 1]{CLCD06} with our Theorem \ref{th:perturbed_theorem}.
The work that is perhaps closest to our work in terms of its convergence analysis is Neely's ``drift-plus-penalty'' algorithm, which, as shown by Huang in  \cite{Hua11}, can be regarded as a randomized version of the dual subgradient method \cite{NB01}. Nonetheless, our analysis is significantly more general since it separates the choice of discrete actions from a specific choice of subgradient update. The latter allows us to construct discrete control policies in a much more flexible manner than previously, which is the second contribution of the paper.  

\subsection{Contribution \textup{(ii)} --- General Control Policies} 

Previous work in \cite{LS04,LS06,ES06,ES07,CLCD06} shows that discrete actions can be obtained as a result of minimizing a linear program. This is very different to our work, where we can have discrete actions in any variable in the optimization and this is not restricted to linear programs.
In our previous work \cite{VL16}, we revisited max-weight policies and provided two different classes of discrete-like updates (Greedy and Frank-Wolfe; see Theorems 1 and 2 in \cite{VL16}). However, those two classes of updates do not provide enough flexibility to design control policies that capture the constraints on the choice of action in many applications; for example, when changes in the action are costly or cannot be made instantaneously. 

In the present work, we show that discrete actions can be selected as any sequence $\{y_k \}$ that stays close, in the appropriate sense, to a sequence $\{x_k\}$ obtained with the subgradient method.  One of the key consequences is that the selection of a discrete action in a time slot can then depend on the action selected in the previous time slot, and so the elements in the sequence of actions can be strongly correlated. 
Another key consequence is that this enables the construction of discrete control policies by only checking that $\| \sum_{i=1}^k x_i - y_i\|_2$ is uniformly bounded. That is, without requiring to re-prove the stability/optimality of the system for every new policy. 
The question of how to select sequence  $\{ y_k \}$ so that $\| \sum_{i=1}^k x_i - y_i\|_2$ remains uniformly bounded is answered in our third contribution.

\subsection{Contribution \textup{(iii)} --- Discrete Control Policies}

In the paper, we propose two new types of discrete control policies: an online policy and another that works with blocks. The policy that works with blocks is more complex but allows us to reshuffle discrete actions, and so we have much more flexibility as to how to select discrete control actions. 

Constraints on the sequence of actions have received relatively little attention in the max-weight literature.   In \cite{CLM12},  \c{C}elik \emph{et al.} show that the original max-weight (myopic) policy fails to stabilize a system when there are reconfiguration delays associated with selecting new actions/configurations, and propose a variable frame size max-weight policy where actions are allocated in frames (\emph{i.e.}, blocks) to minimize the penalties associated with configuration changes. This algorithm is similar to the block-algorithm used in the numerical example in Section \ref{sec:example}, but a notable difference is that in our case the design of the scheduling policy is done simply by ensuring that the $\alpha$-scaled queues stay close to the Lagrange multipliers which yields an approach of much greater generality. 
In \cite{MSY12}, Maguluri \emph{et al.} consider also scheduling in a frame based fashion. However, in their case suboptimal schedules are selected within a frame so that a max-weight schedule can always be selected at the beginning of a frame. Our work differs from \cite{MSY12} fundamentally because we never require to select a max-weight schedule to guarantee optimality. That is, keeping difference $\| \sum_{i=1}^k x_i - y_i\|_2$ uniformly bounded is enough. 
We also note Lin and Shroff's work in \cite{LS06}, where it is observed that the capacity region of the system is reduced as a result of selecting imperfect schedules.\footnote{Imperfect schedules are the result of not being able to select a maximum weight matchings in every time slot.} This is different from our work because we can amortize the complexity of selecting a new schedule/action by selecting the action used in the previous iteration (see Corollary \ref{th:seq_online}). Hence, one does not need to sacrifice resources due to the complexity of selecting a new schedule.

\section{Conclusions}
We have presented a framework that brings the celebrated max-weight features (discrete control actions and myopic scheduling decisions) to the field of convex optimization. In particular, we have shown that these two features can be encompassed in the subgradient method for the Lagrange dual problem by the use of $\delta_k$ and $\epsilon_k$ perturbations, and have established different types of convergence depending on their statistical properties.
One of the appealing features of our approach is that it decouples the selection of a control action from a particular choice of subgradient, and so allows us to design a range of control policies without changing the underlying convex updates. We have proposed two classes of discrete control policies: one that can make discrete actions by looking only at the system's current state, and another which selects actions using blocks. The latter class is useful for handling systems that have constraints on how actions can be selected. We have illustrated this with an example where there are constraints on how packets can be transmitted over a wireless network.

\appendices

\section{Preliminaries}
\label{appendix:preliminaries}

\subsection{Subgradient Convergence}

Let $x_k \in \arg \min_{x \in X} L(x,\lambda_k)$ and observe we can write upper bound
\begin{align*}
h(\lambda_k) - h(\lambda^\star) 
& =   L(x_k, \lambda_k) - L(x^\star,\lambda^\star)  \\
& \le   L(x^\star, \lambda_k) - L(x^\star, \lambda^\star)    \\
& =  (\lambda_k - \lambda^\star)^T g(x^\star)  \\
& \le \| \lambda_k - \lambda^\star \|_2 \| g(x^\star) \|_2 
\end{align*}
and lower bound
\begin{align*}
 h(\lambda_k) - h(\lambda^\star) 
& =   L(x_k, \lambda_k) - L(x^\star,\lambda^\star)   \\
&\ge L(x_k, \lambda_k) - L(x_k, \lambda^\star)   \\
& =  (\lambda_k - \lambda^\star)^T g(x_k)  \\
& \ge -  \| \lambda_k - \lambda^\star \|_2 \| g(x_k) \|_2 .
\end{align*}
If we use the fact that $\| g(x) \|_2 \le \sigma_g$ for all $x \in X$ by Assumption \ref{th:boundedset}, we have
\begin{align}
| h(\lambda_k) - h(\lambda^\star) |  \le \| \lambda_k - \lambda^\star \|_2  \sigma_g. \label{eq:contractionx7az}
\end{align}
From (\ref{eq:contractionx7az}) we can see that if difference $\| \lambda_k - \lambda^\star \|_2$ decreases  then the difference $| h(\lambda_k) - h(\lambda^\star) |$ must eventually also decrease.

\subsection{$\epsilon_k$-subgradients}

In this section, we show that the use of $\epsilon_k$-subgradients is equivalent to approximately minimizing the Lagrangian. Similar results can be found in Bertsekas \cite[pp. 625]{Ber99}, but we include them here for completeness.

Let $x_k \in \arg \min_{x \in X} L(x,\mu_k)$ with $\mu_k = \lambda_k + \epsilon$ for some $\epsilon \in \R^m$, and observe that
$
h(\mu_k) 
 = L(x_k,\mu_k) 
 = L(x_k,\lambda_k+\epsilon) 
 = f(x_k) + (\lambda_k+\epsilon)^T g(x_k) 
 \le f(x_k) + \lambda_k^T g(x_k)  +  \| \epsilon \|_2 \| g(x_k) \|_2
 \le f(x_k) + \lambda_k^T g(x_k)  +  \| \epsilon \|_2 \sigma_g
$
where the last inequality follows since $\sigma_g : = \max_{x \in X} \| g(x) \|_2$. 
Hence, 
\begin{align}
h(\mu_k) - L(x_k,\lambda_k) \le \| \epsilon \|_2 \sigma_g. \label{eq:first_subgrad_low}
\end{align}
We now proceed to show that $| h(\mu_k) - h(\lambda_k) |$ is bounded.  Consider two cases.
Case (i) $h(\mu_k) < h(\lambda_k)$. From the concavity of $h$ we have that 
$
h(\lambda_k) 
 \le h(\mu_k)  + \partial h(\mu_k)^T (\lambda_k - \mu_k) 
 = h(\mu_k)  + \partial h(\lambda_k + \epsilon)^T \epsilon 
 = h(\mu_k)  + g(x_k)^T \epsilon, 
$
and therefore 
\begin{align*}
0 \le h(\lambda_k) - h(\lambda_k + \epsilon)  \le \| \epsilon \|_2 \sigma_g.
\end{align*}
Case (ii) $h(\lambda_k) > h(\mu_k$). Following the same steps than in the first case we obtain
$
h(\mu_k) \le h(\lambda_k,\delta)  - g(x_k)^T \epsilon
$, and therefore
$
0 \le h(\mu_k) - h(\lambda_k)  \le \| g(x_k)  \|_2 \| \epsilon \|_2 \le \| \epsilon \|_2 \sigma_g.
$
Combining both cases yields
$
| h(\lambda_k) - h(\mu_k) | \le  \| \epsilon \|_2 \sigma_g,
$
and using (\ref{eq:first_subgrad_low}) we have
\begin{align*}
0 \le L(x_k,\lambda_k) - h(\lambda_k)  \le  2 \| \epsilon \|_2 \sigma_g := \xi ,
\end{align*}
where the lower bound follows immediately since $h(\lambda_k) \le L(x,\lambda_k)$ for all $x \in X$. 
Hence, the error obtained as a result of selecting an $x_k$ that minimizes $L(x,\mu_k)$ (instead of $L(x,\mu_k)$) is proportional to the difference between $\lambda_k$ and $\mu_k$.

\section{Proofs of Section \ref{sec:main}}
\label{appendix:a}

\subsection{Proof of Lemma \ref{th:perturbed_subgradient_method}}
For any vector $\theta \in \R^m$ we have
\begin{align} 
 \| \lambda_{k+1} - \theta \|_2^2 
& = \| [\lambda_k + \alpha ( g(x_k) + \delta_k) ]^+ - \theta\|_2^2   \notag \\
& \le \| \lambda_k + \alpha ( g(x_k) + \delta_k)  - \theta \|_2^2 \notag \\
& =  \| \lambda_k - \theta\|_2^2 + \alpha^2  \|  g(x_k) + \delta_k \|_2^2  \notag \\
& \qquad + 2\alpha (\lambda_k-\theta)^T  (g(x_k) + \delta_k)  \notag \\
& =  \| \lambda_k - \theta\|_2^2 + \alpha^2  \|  g(x_k) + \delta \|_2^2 + \alpha^2 \| \delta_k - \delta \|_2^2 \notag \\
& \qquad + 2 \alpha^2 (\delta_k - \delta)^T  (g(x_k) + \delta) \notag \\
& \qquad + 2\alpha (\lambda_k-\theta)^T  (g(x_k) + \delta_k)  \label{eq:expansion}
\end{align}
where in the last equation we have used the fact that 
$
\alpha^2  \|  g(x_k) + \delta_k \|_2^2  
 = \alpha^2  \|  g(x_k) + \delta_k - \delta + \delta \|_2^2  
 = \alpha^2     \|  g(x_k) + \delta \|_2^2 + \alpha^2 \| \delta_k - \delta \|_2^2 + 2 \alpha^2 (g(x_k) + \delta)^T  (\delta_k - \delta) .
$
Similarly, observe that 
$
(\lambda_k-\theta)^T  (g(x_k) + \delta_k) 
 = (\lambda_k-\theta)^T  (g(x_k)+\delta + \delta_k - \delta) 
 = (\lambda_k-\theta)^T  (g(x_k) + \delta) + (\lambda_{k} - \theta)^T  (\delta_k - \delta)
$
and since
\begin{align}
&  (\lambda_k-\theta)^T  (g(x_k)+\delta)  \notag \\
& \quad  =(\lambda_k-\theta)^T  (g(x_k)+\delta) + f(x_k) - f(x_k) \notag \\
& \quad  = L (x_k,\lambda_k,\delta) - L(x_k,\theta,\delta) \notag \\
& \quad \le   L (x_k,\lambda_k,\delta) - h(\theta,\delta), \label{eq:saddlemiddle}
\end{align}
we have
\begin{align}
\| \lambda_{k+1} - \theta \|_2^2 
& \le  \| \lambda_k - \theta\|_2^2 + \alpha^2  \|  g(x_k) + \delta \|_2^2 \notag \\
& \quad + \alpha^2 \| \delta_k - \delta \|_2^2 + 2 \alpha^2 (\delta_k - \delta)^T  (g(x_k) + \delta) \notag \\
& \quad + 2\alpha (\lambda_k-\theta)^T  (\delta_k- \delta) \notag \\
& \quad + 2 \alpha ( L (x_k,\lambda_k,\delta) - h(\theta,\delta)) \label{eq:subfullexpansion}
\end{align}
where (\ref{eq:saddlemiddle}) follows from the fact that $h(\theta) = \min_{x \in X} L(x,\theta) \le L(x_k,\theta)$. 
Applying the expansion recursively for $i=1,\dots,k$ we have
\begin{align}
 \| \lambda_{k+1} - \theta \|_2^2 
& \le \| \lambda_1 - \theta\|_2^2  + \alpha^2  \sum_{i=1}^k  \|  g(x_i) + \delta \|_2^2 \notag \\
& \quad +  \alpha^2  \sum_{i=1}^k (\| \delta_i - \delta \|_2^2 + 2  (\delta_i - \delta)^T  (g(x_i) + \delta) ) \notag \\
& \quad + 2\alpha \sum_{i=1}^k  (\lambda_i-\theta)^T  (\delta_i- \delta) \notag \\
& \quad + 2\alpha \sum_{i=1}^k (  L (x_i,\lambda_i,\delta) - h(\theta,\delta) ) 
\end{align}
Next, observe that since 
\begin{align*}
&  L(x_k, \lambda_k,\delta)  \\
&  \quad = L(x_k, \lambda_k,\delta) - L(x_k, \mu_k,\delta) + L(x_k, \mu_k,\delta) \\
 & \quad \le | L(x_k, \lambda_k, \delta) - L(x_k, \mu_k, \delta)| + L(x_k, \mu_k, \delta) \\
& \quad = h(\mu_k, \delta) +  | L(x_k, \lambda_k, \delta) - L(x_k, \mu_k, \delta)| \\
& \quad = h(\mu_k, \delta) +  | (\lambda_k - \mu_k)^T  (g(x_k)+\delta) | \\
& \quad = h(\mu_k, \delta) + |\epsilon_k^T  (g(x_k)+\delta)| \\
& \quad \le h(\mu_k, \delta) +  \| \epsilon_k \|_2 \| g(x_k)+\delta\|_2 \\
& \quad = h(\mu_k, \delta) - h(\lambda_k, \delta)  + h(\lambda_k,\delta) +  \| \epsilon_k \|_2 \|g(x_k) +\delta\|_2 \\
& \quad \le | h(\mu_k, \delta) - h(\lambda_k, \delta)|  + h(\lambda_k, \delta) +  \| \epsilon_k \|_2 \| g(x_k) +\delta\|_2 \\
& \quad \le h(\lambda_k, \delta) + 2 \| \epsilon_k \|_2 \| g(x_k) +\delta\|_2 ,
\end{align*}
we have that
\begin{align}
L(x_k,\lambda_k, \delta) -  h(\lambda_k,\delta) \le  2 \| \epsilon_k \|_2 \| g(x_k) +\delta\|_2 ,
\label{eq:carro}
\end{align}
and therefore
\begin{align}
 \| \lambda_{k+1} - \theta \|_2^2 
& \le  \| \lambda_1 - \theta\|_2^2   + \alpha^2  \sum_{i=1}^k  \|  g(x_i) + \delta \|_2^2 \notag \\
& \quad +  \alpha^2  \sum_{i=1}^k  (  \| \delta_i - \delta \|_2^2 + 2  (\delta_i - \delta)^T  (g(x_i) + \delta) ) \notag  \\
& \quad + 2\alpha \sum_{i=1}^k \left( (\lambda_i-\theta)^T  (\delta_i- \delta) \right) \notag \\
& \quad +  2\alpha \sum_{i=1}^k \left( 2 \| \epsilon_i \|_2 \| g(x_i) +\delta\|_2 \right) \notag \\ 
& \quad +2\alpha  \sum_{i=1}^k (h(\lambda_i,\delta) - h(\theta,\delta))   \label{eq:subfinalexpansion}
\end{align}
Rearranging terms and dividing by $2 \alpha k$ yields the stated result.

\subsection{Proof of Theorem \ref{th:perturbed_theorem}}

Let $\theta = \lambda^\star(\delta)$ in Lemma \ref{th:perturbed_subgradient_method}. From (\ref{eq:subgradient_bound}) and (\ref{eq:carro}) we can write
\begin{align*}
& h(\lambda^\star(\delta),\delta)  \\
& \quad \ge  \frac{1}{k}  \sum_{i=1}^k  h(\lambda_i,\delta)  \\
& \quad \ge \frac{1}{k} \sum_{i=1}^k  \left( L(x_i,\lambda_i,\delta) - 2 \| \epsilon_i \|_2 \| g(x_i) + \delta \|_2 \right) \\
&\quad = \frac{1}{k} \sum_{i=1}^k ( f(x_i) + \lambda_i^T  (g(x_i) + \delta) - 2 \| \epsilon_i \|_2 \| g(x_i) + \delta \|_2) \\
&\quad \ge f(\bar x_k) + \frac{1}{k} \sum_{i=1}^k ( \lambda_i^T  (g(x_i) + \delta) - 2 \| \epsilon_i \|_2 \| g(x_i) + \delta \|_2) , 
\end{align*}
where the last equation follows from the convexity of $f$. 
Rearranging terms
\begin{align}
& f(\bar x_k) - h(\lambda^\star(\delta),\delta) \notag \\
& \quad \le - \frac{1}{k} \sum_{i=1}^k \left( \lambda_i^T  (g(x_i) + \delta) -  2 \| \epsilon_i \|_2 \| g(x_i) + \delta \|_2 \right)
 \label{eq:raw_upper_bound}
\end{align} 
Now, let $\theta = 0$ in (\ref{eq:expansion}) to obtain
\begin{align*}
\| \lambda_{k+1}\|_2^2
& \le  \| \lambda_k \|_2^2 + \alpha^2  \|  g(x_k) + \delta \|_2^2 + \alpha^2 \| \delta_k - \delta \|_2^2 \notag \\
& \quad + 2 \alpha^2 (\delta_k - \delta)^T  (g(x_k) + \delta) + 2\alpha \lambda_k^T  (g(x_k) + \delta_k)
\end{align*}
Using the fact that $\|  g(x_k) + \delta \|_2^2 \le \sigma_g^2$ for all $k$ and applying the latter expansion recursively 
\begin{align}
\| \lambda_{k+1}\|_2^2 & \le  \| \lambda_1 \|_2^2 + \alpha^2  \sigma_g^2  k + \alpha^2 \sum_{i=1}^k \| \delta_i - \delta \|_2^2 \notag \\
& \quad + 2 \alpha^2 \sum_{i=1}^k (\delta_i - \delta)^T  (g(x_i) + \delta) \notag \\
& \quad + 2\alpha \sum_{i=1}^k \lambda_i^T  (g(x_i) + \delta_i)
\end{align}
Rearranging terms, dropping $\| \lambda_{k+1} \|_2$ since it is nonnegative, and dividing by $2\alpha k$ yields
\begin{align*}
- \frac{1}{k} \sum_{i=1}^k \lambda_i^T  (g(x_i) + \delta_i) 
\le  & \frac{\| \lambda_1 \|_2^2}{2 \alpha k} + \frac{\alpha  \sigma_g^2}{2}  + \frac{\alpha}{2k} \sum_{i=1}^k \| \delta_i - \delta \|_2^2 \notag \\
& \quad  + \frac{\alpha}{k} \sum_{i=1}^k (\delta_i - \delta)^T  (g(x_i) + \delta)  
\end{align*}
Combining the last bound with (\ref{eq:raw_upper_bound}), and using the fact that $h(\lambda^\star(\delta),\delta) = f^\star(\delta)$ (by strong duality, \emph{c.f.}, Assumption \ref{th:slater_general}) yields
\begin{align*}
f(\bar x_k) - f^\star(\delta) 
\le  & \frac{\| \lambda_1 \|_2^2}{2 \alpha k} + \frac{\alpha  \sigma_g^2}{2}  + \frac{\alpha}{2k} \sum_{i=1}^k \| \delta_i - \delta \|_2^2  \notag \\ 
& + \frac{\alpha}{k} \sum_{i=1}^k (\delta_i - \delta)^T  (g(x_i) + \delta)  \\
& +  \frac{1}{k} \sum_{i=1}^k 2 \| \epsilon_i \|_2 \| g(x_i) + \delta \|_2
\end{align*}
Taking expectations with respect to $\delta_i$, $i=1,2,\dots,k$ we have $\E (\| \delta_i - \delta \|_2^2) = \sigma_\delta^2$, and $\E((\delta_i - \delta)^T  (g(x_i) + \delta) ) = 0$ since $x_i$ and $\delta_i$ are independent. 
Therefore,
\begin{align}
 \E( f(\bar x_k) - f^\star(\delta)) \le  \frac{\alpha (\sigma^2_g + \sigma_\delta^2)}{2} + \frac{\| \lambda_{1} \|_2^2}{2\alpha k} +  \frac{2}{k} \sum_{i=1}^k  \| \epsilon_i\|_2 \sigma_g \label{eq:upper_bound_raw}
 \end{align}
and we have obtained claim (i). 

We now proceed to lower bound (\ref{eq:upper_bound_raw}). 
Taking expectations with respect to $\delta_i$, $i=1,2,\dots,k$ in Lemma \ref{th:perturbed_subgradient_method}, and using the fact that $\lambda_i$ and $\delta_i$ are independent, so $\E((\lambda_i-\theta)^T (\delta_i-\delta)) = 0$ and $\E((\delta_i-\delta)^T  (g(x_i) + \delta))=0$, we have
\begin{align}
& - \frac{ \| \lambda_1 - \theta\|_2^2 }{2 \alpha k} - \frac{\alpha( \sigma_g^2 + \sigma_\delta^2)}{2} -  \frac{2}{k}\sum_{i=1}^k \| \epsilon_i \|_2 \sigma_g \notag \\
&  \quad  \le   \E \left( \frac{1}{k} \sum_{i=1}^k h( \lambda_i,\delta)  - h(\theta,\delta) \right).
 \label{eq:big_sum_exp}
\end{align}
Next, by the convexity of $-h(\cdot,\delta)$ we can write
\begin{align*}
\frac{1}{k}  \sum_{i=1}^k  \E  ( h(\lambda_i,\delta) )
=  \E \left( \frac{1}{k}  \sum_{i=1}^k h(\lambda_i,\delta) \right) 
\le \E(h(\bar \lambda_k),\delta)
\end{align*}
and letting $\theta = \lambda^\star (\delta)$ 
\begin{align}
&   - \frac{\| \lambda_1 - \lambda^\star(\delta)\|_2^2}{2\alpha k}  - \frac{\alpha ( \sigma_g^2 + \sigma_\delta^2) }{2} -  \frac{2}{k} \sum_{i=1}^k \| \epsilon_i \|_2 \sigma_g \notag  \\
& \quad \le  \E \left( h(\bar \lambda_k,\delta) - h(\lambda^\star(\delta),\delta) \right)  \le 0, \label{eq:avg_h_bound}
\end{align}
where the upper bound follows from the fact that $h(\lambda^\star(\delta), \delta) = \sup_{\lambda \succeq 0} h(\lambda,\delta)$. Next, from the saddle point property of the Lagrangian 
\begin{align*}
\E(h(\bar \lambda_k),\delta)
& \stackrel{\text{(a)}}{\le} \E(L(\E(\bar x_k), \bar \lambda_k, \delta)) \\
& = \E ( f(\E(\bar x_k)) + \E (\bar \lambda_k)^T  (g(\E(\bar x_k)) + \delta)) ) \\
& \stackrel{\text{(b)}}{\le} \E(f(\bar x_k)) + \E (\bar \lambda_k)^T  \E(g(\bar x_k) + \delta),
\end{align*} 
where the expectation on $\bar x_k$ in the RHS of (a) is taken with respect to $\delta_{i}$, $i=1,\dots,k$; and (b) follows from the convexity of $f$ and $g$. Therefore,
\begin{align}
&   - \frac{\| \lambda_1 -  \lambda^\star(\delta)\|_2^2}{2\alpha k}  - \frac{\alpha ( \sigma_g^2 + \sigma_\delta^2) }{2} -  \frac{2}{k} \sum_{i=1}^k \| \epsilon_i \|_2 \sigma_g  \notag \\
& -  \E(\bar \lambda_k)^T  \E(g(\bar x_k) + \delta)   \le \E( f(\bar x_k) - f^\star(\delta)). \label{eq:raw_lower_bound}
\end{align}
We need to show that $\E(\bar \lambda_k)^T  \E (g(\bar x_k) + \delta)$ is upper bounded. 
Observe first that for any sequence $\{ x_k \}$ from $X$ we can write
\begin{align*}
\lambda_{k+1} = [\lambda_{k} + \alpha (g(x_k) + \delta_k) ]^+  \succeq  \lambda_k + \alpha (g(x_k)+ \delta_k ),
\end{align*}
and applying the latter recursively we have that
\begin{align*}
\lambda_{k+1}  \succeq \lambda_1 + \alpha \sum_{i=1}^k  (g(x_i) + \delta_i).
\end{align*}
Dropping $\lambda_1$ since it is nonnegative, dividing by $\alpha k$, and using the convexity of $g$ follows that
\begin{align*}
g(\bar x_k) + \frac{1}{k}\sum_{i=1}^k \delta_i \preceq \frac{\lambda_{k+1}}{\alpha k}, 
\end{align*}
and taking expectations with respect to $\delta_i$, $i=1,\dots,k$
\begin{align}
\E( g(\bar x_k) + \delta)) \preceq \frac{\E(\lambda_{k+1})}{\alpha k}.
\label{eq:constraint_violation_bound_proof}
\end{align}
Multiplying both sides by $\E(\bar  \lambda_k)$ (where the expectation is with respect to $\delta_i$, $i=1,\dots,k$) and using Cauchy-Schwarz  
\begin{align}
\E(\bar  \lambda_k)^T  \E( g(\bar x_k) + \delta) 
& \le \frac{\E(\bar  \lambda_k)^T  \E(\lambda_{k+1})}{\alpha k} \notag \\
&  \le \frac{ \| \E(\bar  \lambda_k) \|_2 \| \E(\lambda_{k+1}) \|_2}{\alpha k} . \notag
\end{align}

We proceed to show that $\| \E(\bar \lambda_k) \|_2$ is bounded using Lemma 6 in \cite{VL16}. This lemma says that for any $\chi \ge 0$ then $\mathcal Q_\chi :=\{ \lambda \succeq 0 \mid  h(\lambda) \ge h(\lambda^\star(\delta)) - \chi\}$ is a bounded set. Further, for any $\lambda \in \mathcal Q_\chi$ we have that $\| \lambda \|_2 \le \frac{1}{\upsilon}(f(\hat x) - h(\lambda^\star(\delta)) + \chi)$ where $\hat x$ is a Slater point, and $\upsilon := \min_{j \in \{1,\dots,m\}} -g_j(\hat x)$ a constant that does not depend on $\chi$. Note that $\upsilon > 0$.
Now, observe that since $\E( h( \bar \lambda_k,\delta) ) \le h(\E(\bar \lambda_k),\delta)$, from (\ref{eq:avg_h_bound}) we can write
\begin{align}
&  - \frac{ \| \lambda_1 -  \lambda^\star(\delta)\|_2^2 }{2 \alpha k} - \frac{\alpha ( \sigma_g^2 + \sigma_\delta^2)}{2} -  \frac{2}{k} \sum_{i=1}^k \| \epsilon_i \|_2  \sigma_g \notag \\
& \quad \le  h( \E (\bar \lambda_k),\delta) - h(\lambda^\star(\delta),\delta)  \le 0.
\label{eq:avg_rv_cvg}
\end{align}
Hence, if we identify $-\chi$ with the LHS of (\ref{eq:avg_rv_cvg}) we obtain that $\|\E (\bar \lambda_k)\|_2$ is bounded. That is, 
\begin{align}
\| \E (\bar \lambda_k) \|_2 
& \le \frac{1}{\upsilon} \Bigg(  f(\hat x) - h(\lambda^\star(\delta)) +  \frac{ \| \lambda_1 -  \lambda^\star(\delta)\|_2^2 }{2 \alpha k} \notag \\ 
& \qquad   \qquad + \frac{\alpha ( \sigma_g^2 + \sigma_\delta^2)}{2} +  \frac{2}{k} \sum_{i=1}^k \| \epsilon_i \|_2  \sigma_g \Bigg) \notag \\
& \le \frac{1}{\upsilon} \left(  f(\hat x) - h(\lambda^\star(\delta)) + \frac{\Omega}{\alpha k} \right) \label{eq:multiplier_bound_xx}.
\end{align}
We continue by giving a bound on $\| \E(\lambda_{k+1}) \|_2$. Taking expectations in (\ref{eq:subfinalexpansion}) with respect to $\delta_i, \ i=1,\dots,k$, letting $\theta = \lambda^\star(\delta)$, $\| g(x_k) + \delta \|_2 \le \sigma_g$, and using the fact that $\lambda_k$ and $\delta_k$ are independent for all $k$, we have
\begin{align}
& \E(\| \lambda_{k+1} - \lambda^\star(\delta) \|_2^2 ) \notag \\
& \quad \le  \| \lambda_1 - \lambda^\star(\delta)\|_2^2  + \alpha^2  (\sigma_g^2 + \sigma_\delta^2)  k  + 2 \alpha \sum_{i=1}^k \| \epsilon_i \|_2 \sigma_g  \notag \\
& \qquad + 2\alpha  \sum_{i=1}^k (h(\lambda_i,\delta) - h(\lambda^\star(\delta),\delta))  \notag 
\end{align}
Next, observe that since $h(\lambda_i,\delta) - h(\lambda^\star(\delta),\delta) \le 0$ for all $i=1,\dots,k$ we can write
$
\E(\| \lambda_{k+1} - \lambda^\star(\delta) \|_2^2 )
 \le  \| \lambda_1 - \lambda^\star(\delta)\|_2^2  
+ \alpha^2  (\sigma_g^2 + \sigma_\delta^2)  k  + 2 \alpha \sum_{i=1}^k \| \epsilon_i \|_2 \sigma_g   
$
and by using the convexity of $\| \cdot\|_2^2$,  
$
\| \E(\lambda_{k+1}) - \lambda^\star(\delta) \|_2^2
 \le  \| \lambda_1 - \lambda^\star(\delta)\|_2^2  
+ \alpha^2  (\sigma_g^2 + \sigma_\delta^2)  k  + 2 \alpha \sum_{i=1}^k \| \epsilon_i \|_2 \sigma_g. 
$
That is, $\E(\lambda_{k+1})$ is within a ball around $\lambda^\star(\delta)$. 
Next, since $\| \lambda^\star(\delta) \|_2^2$ is bounded we can write 
$
\| \E(\lambda_{k+1})  \|_2^2
 \le \| \lambda^\star(\delta) \|_2^2 +  \| \lambda_1 - \lambda^\star(\delta)\|_2^2  
+ \alpha^2  (\sigma_g^2 + \sigma_\delta^2)  k  
+ 2 \alpha \sum_{i=1}^k \| \epsilon_i \|_2 \sigma_g  
$
and by taking square roots in both sides and using the concavity of $x^a$ for $x \in \R_{++}$ with $0 < a < 1$ (see \cite[pp. 71]{BV04}) we have 
\begin{align}
\| \E(\lambda_{k+1})  \|_2
& \le \| \lambda^\star(\delta) \|_2 +  \Big( \| \lambda_1 - \lambda^\star(\delta)\|_2^2 \notag \\
& \qquad + \alpha^2  ( \sigma_g^2 + \sigma_\delta^2) k  + {2 \alpha \sum_{i=1}^k \| \epsilon_i \|_2 \sigma_g} \Big)^{-1/2} \notag \\
& =  \| \lambda^\star(\delta) \|_2  + \sqrt \Omega 
\label{eq:lucky_bound}
\end{align}
where the last inequality follows from the concavity of the square root and the fact that all the terms are nonnegative.
Hence, 
\begin{align*}
& \E(\bar  \lambda_k)^T  \E( g(\bar x_k) + \delta) \le \| \E(\bar  \lambda_k) \|_2  \left( \frac{ \| \lambda^\star(\delta) \|_2}{\alpha k}  + \frac{\sqrt \Omega}{\alpha k} \right) \end{align*}
and so we can use (\ref{eq:raw_lower_bound}) to lower bound (\ref{eq:upper_bound_raw}), and obtain the bound claimed in (ii). 

Claim (iii) follows from (\ref{eq:constraint_violation_bound_proof}) and (\ref{eq:lucky_bound}).
Claim (iv) follows from (\ref{eq:multiplier_bound_xx}).

\section{Proof of Section \ref{sec:actions}}
\label{appendix:c}
\subsection{Proof of Lemma \ref{th:blocksequence}}
To start, let $V = |Y|$ and note we always have $\1^T  (z + \delta) = V$ since $z$ is the sum of $V$ elements from $U$,  $u^T  \1 = 1$ and $\delta^T  \1 = 0$ for all $u\in U$, $\delta \in D$. Further, $(z+\delta) \succeq - \1$ since $\delta \succeq -\1$ and $z \succeq 0$. 
Now, let $r:=(z + \delta)$ and define 
\begin{align*}
a = -[-r]^+, \quad b = \lfloor r - a \rfloor, \quad   c = r - a-b,
\end{align*}
where the floor in $b$ is taken element-wise. 
That is,  $a \in [-1,0]^V$, $b \in \{0,1,\dots,V\}^V$, $c \in [0,1)^V$. 
For example, if $r = [2.2,-0.2]^T $ then $a = [0,-0.2]^T $, $b = [2,0]^T $ and $c = [0.2,0]^T $.
Observe, 
\begin{align*}
\1^T  r = \1^T  (a + b + c) = V,
\end{align*}
and since $b$ is integer valued $\1^T  b \in \Z_+$, which implies $\1^T  (a+c) \in \Z_+$. 
Next,  let $\1^T  b = V-\1^T (a+c) := V'$, and observe $b$ can be written as the sum of $V'$ elements from $E$, \emph{i.e.}, 
\begin{align*}
b = \sum_{i=1}^{V'} e_i.
\end{align*} 
Next, since $-\1 \preceq a+c \prec \1$ and $\1^T  (a+c) = V - V':=V''$, there must exist \emph{at least} $V''$ elements in vector $(a+c)$ that are nonnegative. 
If we select $V''$  elements from $E$ that match the nonnegative components of vector $(a+c)$ we can construct a subsequence $\{ e_i \}_{i=1}^{V''}$ such that 
\begin{align*}
- \1 \preceq (a+c) - \sum_{i=1}^{V''} e_i \prec \1 .
\end{align*}
Finally, letting $z' = \sum_{i=1}^{V'} e_i + \sum_{i=1}^{V''} e_i $ yields the result. 

\subsection{Proof of Lemma \ref{th:magic_update_lemma}}

First of all, recall from the proof of Lemma \ref{th:blocksequence} that 
$\1^T  (\delta  + z - z') = 0$, $\delta + z \succeq -\1$ and that $V:= \1^T (z + \delta)$ where 
$V : = |Y|$. Next, define 
$r_i := \delta + z - \sum_{\kappa=1}^{i-1} e_i$, $i=1,2,\dots$ and note that update $
e_i \in \arg \min_{e \in E} 
 \| r_{i} - e \|_\infty
$
decreases the largest component of vector $r_{i}$, \emph{i.e.}, in each iteration a component of vector $r_i$ decreases by $1$, and therefore $\1^T  r_i = V - i + 1$ with $i = 1,\dots, |Y|+1$.

For the lower bound observe that if $r_{i+1}(j) < -1$ for some $j = 1,\dots, |Y|$ we must have that $r_{i} \prec 0$ since the update $e_i \in \arg \min_{e \in E} 
 \| r_{i} - e \|_\infty$ selects to decrease the largest component of vector $r_{i}$. 
However, since $\1^T  r_{i} \ge  0$ for all $i= 1,\dots,|Y|+1$ we have that vector $r_{|Y|}$ has at least one component that is nonnegative. Therefore, $r_{|Y|+1} \succeq -\1$ and $\delta + z - z' \succeq - \1$.
For the upper bound define $a_{i} = -[-r_{i}]^+$,  $b_{i} = \lfloor r_{i} - a_i \rfloor$, $ c_{i} = r_{i} - a_{i}-b_{i}$, $i=1,\dots,|Y|+1$ and note that $-\1 \preceq a_{i} \preceq 0$ and $ 0 \preceq c_{i} \prec \1$ for all $i=1,\dots,|Y| + 1$, and that $\1^T  b_{i}$ decreases by $1$ in each iteration if $\1^T  b_i \ge 1$. Hence, $b_{|Y|+1} = 0$ and therefore $-\1 \preceq r_{|Y|+1} = a_{|Y|+1} + c_{|Y|+1} = \delta + z - z' \preceq \1$ and we are done.

\subsection{Proof of Theorem \ref{th:online_actions}}

We begin by noting that since $\1^T u_k=1=\1^T  e_k$ then $\1^T  s_k=0$ for all $k=1,2,\dots$.  Also note that since $u_k \in U$ all elements of $u_k$ are nonnegative and at least one element must be non-zero since $\1^T  u_k=1$.

We now proceed by induction to show that there always exists a choice of $e_{k+1}$ such that $s_k \succeq -\1$, $k=1,2,\dots$.  
When $k=1$ let element $u_1(j)$ be positive (as already noted, at least one such element exists).  Selecting $e_1 = v_j$ then it follows that $-1< u_1(j) - e_1(j)\le 0$ and so $-1\prec u_1- e_1 \prec 1$.  That is, ${s}_1\succeq -\1$.   
Suppose now that ${s}_k\succ -\1$.  We need to show that ${s}_{k+1}\succeq -\1$.   Now ${s}_{k+1}= {s}_k+ {u}_{k+1}- {e}_{k+1}$.   Since ${s}_k \succeq -\1$, ${s}_k (j) \ge -1$ $\forall j=1,\dots, |Y|$.   Also, $\1^T  s_k = 0$, so either all elements are 0 or at least one element is positive.   If they are all zero then we are done (we are back to the $k=1$ case).    
Otherwise, since all elements of ${u}_{k+1}$ are nonnegative then at least one element of ${s}_k+ {u}_{k+1}$ is positive.  Let element $s_k(j)+ u_{k+1}(j)$ be the largest positive element of $s_k+u_{k+1}$.   Selecting $e_{k+1}= {v}_j$ then it follows that $s_k(j)+u_{k+1}(j)-e_{k+1}(j)\ge -1$.   
That is, $s_{k+1}\succeq -\1$.

We now  show that $s_k$ is upper bounded. Recall ${e}_{k+1}$ can always be selected such that ${s}_k \succeq - \1$, and also $\1^T  {s}_k=0$.  Since $\1^T  {s}_k=0$ either $s_k$ is zero or at least one element is positive. Since $\vv{s}_k\succeq -\1$ and at most $|Y|-1$ elements are negative, then the sum over the negative elements is lower bounded by $-(|Y|-1)$. Since $\1^T  {s}_k=0$ it follows that the sum over the positive elements must be upper bounded by $|Y|-1$. 
Hence, $\| {s}_k\|_\infty \le (|Y|-1)$.

\subsection{Proof of Corollary \ref{th:seq_online}}

Since $s_k$ has at least one component that is nonnegative, and update (\ref{eq:ekupdate}) selects the largest component of vector $s_k$ when $k \in \mathcal T$, we have that a component of vector $s_k$ can decrease at most by $\bar \tau$ in an interval $\{ \tau_j- \tau_{j+1}\}$ for all $j=1,2,\dots$. Hence, $s_k \succeq - \bar \tau \1$ for all $k$. 
Next, since $s_k^T  \1 = 0$ for all $k$ and the sum over the negative components is at most $-\bar \tau (|Y| - 1)$, we have that $ s_k \preceq \bar \tau  (|Y|-1)  \1$. The rest of the proof follows as in  Theorem \ref{th:online_actions}.

\subsection{Proof of Theorem \ref{th:gen_actions}}

Recall that since $\1^T u_k=1=\1^T  e_k$ then $\1^T  s_k=0$ for all $k=1,2,\dots$, and therefore $s_k$ is either $0$ or at least one of its components is strictly positive. 
Next, observe that since $\gamma_k = - \min_{j \in \{1,\dots,|Y| \} } s_k(j)$ we have that $\max_{j \in \{1,\dots,|Y| \} } s_k(j) \le \gamma_k (|Y|-1)$, which corresponds to the case where $|Y|-1$ components of vector $s_k$ are equal to $\gamma_k$. The rest of the proof continues as in the proof of Theorem \ref{th:online_actions}.

\balance
\bibliographystyle{IEEEtran}
\bibliography{references}
%

%

\end{document}